\documentclass[12pt]{elsarticle}
\pdfoutput=1

 \makeatletter
\def\bs{\expandafter\@gobble\string\\}
\def\lb{\expandafter\@gobble\string\{}
\def\rb{\expandafter\@gobble\string\}}
\def\@pdfauthor{C.V.Radhakrishnan}
\def\@pdftitle{elsarticle.cls -- A documentation}
\def\@pdfsubject{Document formatting with elsarticle.cls}
\def\@pdfkeywords{LaTeX, Elsevier Ltd, document class}

\DeclareRobustCommand{\LaTeX}{L\kern-.26em%
        {\sbox\z@ T%
         \vbox to\ht\z@{\hbox{\check@mathfonts
           \fontsize\sf@size\z@
           \math@fontsfalse\selectfont
          A\,}%
         \vss}%
        }%
     \kern-.15em%
    \TeX}
\makeatother

\usepackage{graphicx}
\usepackage[T1]{fontenc}

\usepackage{amsmath}
\usepackage{amssymb,amsthm}

\usepackage{bm}
\usepackage{xcolor}
\usepackage{subfigure}

\usepackage{float}

\begin{document}

\title{The role of the initial distribution in population survival within a bounded habitat}

\author{Rafael de la Rosa\corref{cor1}\fnref{fn1}}
\ead{rafael.delarosa@uca.es}
\author{Elena Medina\fnref{fn1}}
\ead{elena.medina@uca.es}
\cortext[cor1]{Corresponding author}
\fntext[fn1]{Departamento de Matem\'aticas,
                        Facultad de Ciencias,
                        Universidad de C\'adiz,
                        11510 Puerto Real, C\'adiz, Spain}

\begin{abstract}

In this paper, we analyze the role of initial conditions in population persistence. Specifically, we consider the reaction-diffusion equation $u_t\,=\,D\,(u^{\nu-1}\,u_x)_x\,+\,a\,u^{\mu}$, with $\mu,\nu>0$, accompanied by hostile boundary conditions and examine two families of one-parametric initial distributions, including homogeneous distributions. The model was previously studied by Colombo and Anteneodo (2018). They determined appropriate habitat sizes $l$ for the survival of a population, whose individuals are initially placed homogeneously within the full habitat domain with a total initial population $n_0$.
We show that the survival condition can be naturally formulated in terms of the parameter $Q:=\frac{a}{D}l^{-\mu+\nu+2}n_0^{\mu-\nu}$. Indeed, there exists a critical value $Q_c$ determined by $\mu$, $\nu$ and the initial distribution parameter such that the survival condition can always be written as $Q\geq Q_c$. Notably, from this point of view, one can derive a condition for $Q$ that holds universally for our model under conditional persistence ($\mu\geq\nu$). It applies, in particular, to the case $\mu=\nu+2$, which was not addressed in the previously mentioned work. Nevertheless, in this case $Q=\frac{a}{D}n_0^2$, therefore survival depends solely on the total population, not on the habitat size. We apply a finite-difference scheme to estimate $Q_c$. Conversely, given a population whose evolution is determined by $\mu$, $\nu$, $l$, $n_0$, and the growth and diffusion coefficients $a$ and $D$ (and consequently the value of $Q$) we use the numerical algorithm to estimate the initial distribution to ensure population survival.
\end{abstract}


\date{today}
\maketitle
\section{Introduction\label{sec:intro}}
The analysis of the evolution of a population is a crucial topic in biology, ecology, medicine, agriculture,.... One of the main concerns is the survival of the species, and more specifically, the conditions that the habitat has to fulfill for the persistence of a population. Since the pioneering Skellam's work~\cite{skellam1951random}, the effect of the size of the habitat on the fate of the population has been thoroughly investigated~\cite{berti2015extinction,cantrell2002habitat,colombo2018nonlinear,colombo2016population,dornelas2024movement,dos2020critical,holmes1994partial,latorre1998spatial,
li2022wave,lin2004localization,ludwig1979spatial,maciel2018critical,neicu2000extinction,
perry2005experimental,xu2020population,zhou2017discrete}. The problem is usually posed in terms of diffusion-reaction equations, however, convective terms are usually included (see, for example,~\cite{holmes1994partial} where the effects of the wind or water currents are taken into account,~\cite{neicu2000extinction} where the extincion velocity of bacteria colonies under forced convection is analyzed,~\cite{lin2004localization} and more recently~\cite{dornelas2024movement}). In these models the effect of the boundary can be completely hostile, i.e., individuals die if they reach the habitat boundary, and consequently, Dirichlet homogeneous boundary conditions are imposed. Nevertheless, different degrees to which the region outside the patch is letal are also considered. Some of the papers devoted to populations in bounded habitats with only partially hostile boundaries are:~\cite{ludwig1979spatial}
where the width of the partially hostile sourronding region to avoid an outbreak state is studied, ~\cite{jin2023enhancing} where the optimal location of a protective zone is investigated, ~\cite{cantrell2002habitat} for a predator-prey model,~\cite{maciel2018critical} for a two-sex population model, or~\cite{xu2020population} for a competition model. Other alternative models are discussed in~\cite{colombo2016population}
where the patch size is allowed to depend on time,~\cite{dos2020critical} where the diffusion coefficient is, in general, space dependent, or~\cite{dornelas2024movement} where the models include the existence of preferred regions in the habitat.

Most of these models consider a linear diffusion term, and conclude the existence of a critical patch size such that the population survives for bigger habitats sizes whereas it gets extincted for smaller habitat sizes. Experimental confirmation of these results is provided in~\cite{perry2005experimental} for an E-coli population in a quasi-one dimensional habitat.
It is also worth to mention~\cite{holmes1994partial} where two spatial  dimensional models are studied and it is found that the smallest critical patch areas correspond to circular habitats while the critical area increases as the habitat geometry deviates from circularity. Another remarkable paper in which not only the habitat size but also the region geometry are found to determine the survival conditions is~\cite{tam2022effect}, where a Fisher-Stefan model with a moving boundary is analyzed.

Similar results are obtained by using discrete models~\cite{berti2015extinction} or models in terms of integro-difference equations~\cite{latorre1998spatial,li2022wave,zhou2017discrete}.

In general, nonlinear models described by partial differential equations (PDEs) cannot be solved by using purely analytical methods because their complexity, and exact solutions are only known for certain particular cases. Furthermore, the utility of these solutions is further limited to models involving boundary conditions which must be verified, even though these boundary conditions are defined by equations.

In such situations, one can employ numerical or analytical approximation methods which can yield crucial insights into the behaviour of the solution at critical values of the parameters that the PDE involves. Nevertheless, analytical approximations methods are typically more challenging to apply than numerical methods.

Among the approximation tecnhiques available for solving PDEs, finite-difference and finite element methods stand out for their versatility and applicability. In this paper, we use a finite-difference method. Finite-difference methods for PDEs approximate partial derivatives at a point by difference quotients in a small neighborhood of the point. These methods provide in most cases solutions that are either as accuarate as the data permit or sufficiently precise for the aim for which the solutions are obtained. Solutions computed using finite-difference methods are as suitable as those obtained from an analytical formula.

In~\cite{colombo2018nonlinear}, Colombo and Anteneodo consider the determination of appropriate habitat sizes for the persistence of a population which evolves according to a model with nonlinear diffusion term. They assume that the habitat can be considered one dimensional and the concentration of individuals obeys to the reaction-diffusion equation
\begin{equation}\label{eq1}
u_t\,=\,D\,(u^{\nu-1}\,u_x)_x\,+\,a\,u^{\mu},
\end{equation}
where both diffusion and growth coefficients are not constant but depend on the concentration through a power of the density with a specific exponent~\cite{piva2021interplay}.
The population is supposed to be confined in a habitat of size $l$, such that outside the habitat there is not appropriate life conditions for 
the individuals so that if they reach the boundary, they die. Thus, the boundary conditions can be written as
\begin{equation}\label{boun}
u\left(-\frac{l}{2},t\right)\,=\,u\left(\frac{l}{2},t\right)\,=\,0,
\end{equation}
and the initial distribution is assumed to be homogeneous with a total initial population $n_0$, i.e. 
\begin{equation}\label{ini}
u(x,0)\,=\,\frac{n_0}{l}.
\end{equation}
The results in~\cite{colombo2018nonlinear} can be sumarized as
\begin{enumerate}
\item For $\mu\,<\,\nu$ the population survives independently of the habitat size $l$ and the total initial population $n_0$. Moreover, it reaches a steady state.
\item For $\mu\,=\,\nu$, there exits a critical habitat size $l_c$ which does not depend on the total initial population. The population gets extinted for smaller habitat sizes whereas it grows unboundedly for larger habitat sizes.
\item For $\mu\,\in\,(\nu,\nu\,+\,2)$, the survival condition is similar as in the previous case but the critical habitat sizes does depend on the total initial population through a factor of the
form $n_0^{\frac{\nu-\mu}{2+\nu-\mu}}$.
\item For $\mu\,>\,\nu\,+\,2$, there also exists a critical habitat size depending on the total initial population through a factor $n_0^{\frac{\nu-\mu}{2+\nu-\mu}}$. However, in this case the population gets extinted for larger habitat sizes while it survives (with unbounded growth) for smaller habitat sizes.
\end{enumerate}

The purpose of this work is to analyze the role of the initial conditions for the model given by the reaction-diffusion equation~(\ref{eq1}) with boundary conditions (\ref{boun}) in the case of the conditional persistence $\mu\,\geq\,\nu$. Indeed, this can be a relevant problem when a habitat with a given size has to be repopulated by a certain endangered species and only a fixed amount of individuals are available. 

\begin{figure}[H]
	\centering
        \includegraphics[width=10cm]{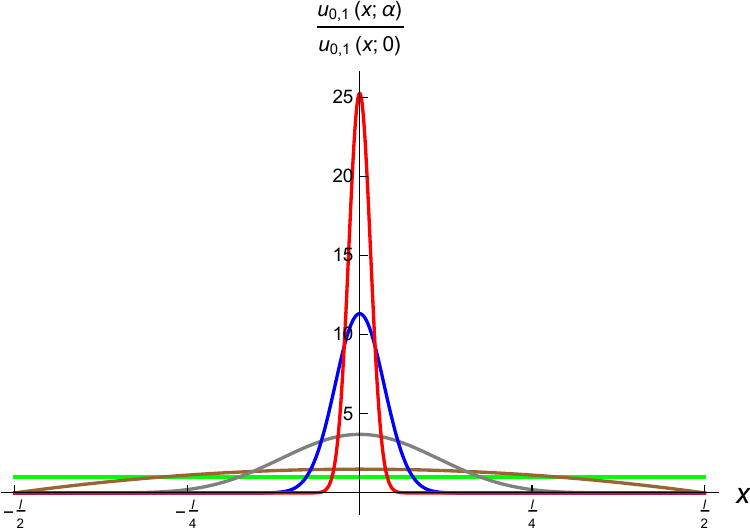}
	\caption{Initial distributions~(\ref{u01}) for the values $\alpha=0$ (green, homogeneous distribution), $\alpha=1$ (brown), $\alpha=10$ (gray), $\alpha=100$ (blue) and $\alpha=500$ (red).
	\label{fig:u01}}
\end{figure}

\begin{figure}[H]
	\centering
        \includegraphics[width=10cm]{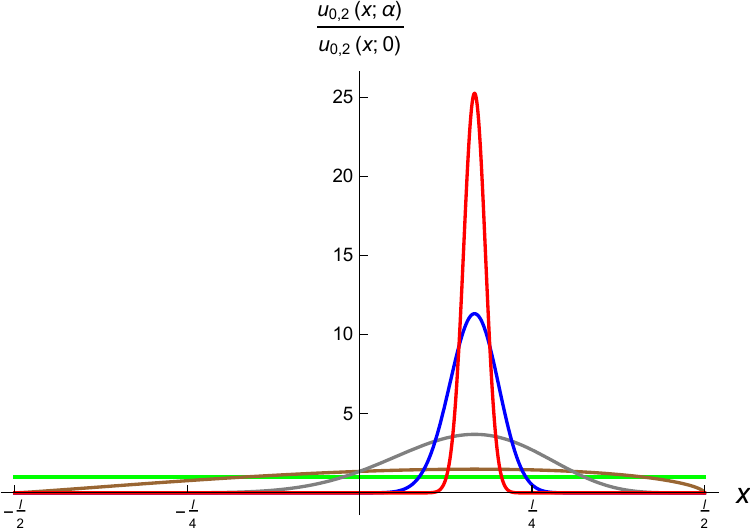}
	\caption{Initial distributions~(\ref{u02}) for the values $\alpha=0$ (green, homogeneous distribution), $\alpha=1$ (brown), $\alpha=10$ (gray), $\alpha=100$ (blue) and $\alpha=500$ (red).
	\label{fig:u02}}
\end{figure}

The starting point is the problem 
\begin{equation}\label{prob1}\left\{\begin{array}{ll}
u_t\,=\,D\,(u^{\nu-1}\,u_x)_x\,+\,a\,u^{\mu}, &  x\,\in\,\left[-\frac{l}{2},\frac{l}{2}\right], \;  t\,>\,0,\\  \\
u\left(-\frac{l}{2},t\right)\,=\,u\left(\frac{l}{2},t\right)\,=\,0,     &t\,>\,0,\\  \\
u(x,0)\,=\,u_0(x),  &  x\,\in\,\left[-\frac{l}{2},\frac{l}{2}\right],
\end{array}\right.\end{equation}
where the total initial population is given by
\begin{equation}\label{n0}
\int_{-\frac{l}{2}}^{\frac{l}{2}}u_0(x)\,dx\,=\,n_0.
\end{equation}
First, using nondimensional variables we find that the problem is naturally and conveniently formulated in terms of the parameter
\begin{equation}\label{Q}
Q:=\frac{a}{D}l^{-\mu+\nu+2}n_0^{\mu-\nu}.
\end{equation}
More particularly, there exists a critical value of $Q$ which we denote by $Q_c$ such that the survival condition can always be written as $Q\,\geq\,Q_c$, and these critical values depends 
on the exponents $\mu$ and $\nu$ and on the initial distribution $u_0(x)$. In order to analyze these dependencies we consider families of one-parametric initial distributions $u_0(x;\alpha)$ corresponding to a habitat size $l$, a total population $n_0$ and a certain parameter $\alpha$ such that
\begin{itemize}
\item $\alpha$ takes non-negative values.
\item $u_0(x;0)=\frac{n_0}{l}$, i.e., for $\alpha=0$ the initial distribution is the homogeneous initial distribution.
\item As $\alpha$ grows the initial distribution is more and more concentrated in a small neighborhood of 
a point of the habitat size (see Figures~\ref{fig:u01} and \ref{fig:u02}).
\end{itemize}

Then, we can apply a numerical algorithm to analyze how the critical value depends on the exponents $\mu$ and $\nu$, and on the parameter $\alpha$ in the initial condition.
Reciprocally, given the parameters in the population, i.e., the exponents $\mu$ and $\nu$, the habitat size $l$, the total population $n_0$ and the growth and diffusion coefficients $a$ and $D$ (and consequently the value of $Q$) we use the numerical algorithm to estimate
the value of  the parameter $\alpha$ or more qualitative ``how concentrated the initial population should be'' to ensure survival.

The first family of initial conditions we consider is given by
\begin{equation}\label{u01}
u_{0,1}(x;\alpha)\,=\,\frac{1}{B(1+ \alpha,1+\alpha)}\left(\frac{1}{4}-\frac{x^2}{l^2}\right)^{\alpha}\frac{n_0}{l}, \quad \alpha\geq0,
\end{equation}
where $B(a,b)$ denotes the Beta function. Note that $u_{0,1}(x;0)\,=\,n_0/l$ corresponds to the homogeneous distribution and as $\alpha$ increases, the initial distributions tends to concentrate near the origin. Furthermore, we consider a second family of initial conditions given by
\begin{equation}\label{u02}
u_{0,2}(x;\alpha)=\frac{1}{B(1+\gamma(\alpha),1+2\gamma(\alpha))}\left(\left(\frac{1}{2}-\frac{x}{l}\right)\left(\frac{1}{2}+\frac{x}{l}\right)^2\right)^{\gamma(\alpha)}\frac{n_0}{l}, \quad \alpha\geq0,
\end{equation}
where $\gamma(\alpha)$ is determined as the unique solution of the equation
\begin{equation}\label{gammaalpha}
4^{\alpha}B(1+ \alpha,1+\alpha)\,=\,\left(\frac{27}{4}\right)^{\gamma}B(1+\gamma,1+2\gamma).
\end{equation}
Note that, according to~(\ref{gammaalpha}) $\gamma(0)=0$, so that $u_{0,2}(x,0)=n_0/l$. Moreover, $u_{0,2}(\cdot;\alpha)$ reaches its maximum at $x=l/6$, and consequently,
\begin{equation}\label{u02max}
u_{0,2,max}(\alpha)\,=\,u_{0,2}\left(\frac{l}{6};\alpha\right) \,=\, \left(\frac{4}{27}\right)^{\gamma(\alpha)}\frac{1}{B(1+\gamma(\alpha),1+2\gamma(\alpha))}\frac{n_0}{l}.
\end{equation}
On the other hand, the first family of initial distributions is symmetric with respect to the midpoint of the habitat, $x=0$, so that
\begin{equation}\label{u01max}
u_{0,1,max}(\alpha)\,=\,u_{0,1}(0;\alpha) \,=\, \left(\frac{1}{4}\right)^{\alpha}\frac{1}{B(1+\alpha,1+\alpha)}\frac{n_0}{l}.
\end{equation}
Thus, for the choice of $\gamma(\alpha)$ given by equation~(\ref{gammaalpha}), both families of initial distributions satisfy 
\begin{equation}\label{eqmax}
u_{0,1,max}(\alpha)\,=\,u_{0,2,max}(\alpha).
\end{equation}
We plot the ratios of the initial distributions $u_{0,1}(\cdot;\alpha)$ (respectively $u_{0,2}(\cdot;\alpha)$) to the homogeneous initial distribution $n_0/l$ in Figure~\ref{fig:u01} (respectively Figure~\ref{fig:u02}) for the values $\alpha=0$ (green, homogeneous distribution), $\alpha=1$ (brown), $\alpha=10$ (gray), $\alpha=100$ (blue) and $\alpha=500$ (red).

This paper is organized as follows. In Section \ref{sec:model}, we present a model given by a reaction-diffusion equation for the evolution of the population density in a bounded habitat. We show that using nondimensional variables the survival condition can always be naturally expressed as $Q \geq Q_c$, with $Q:=\frac{a}{D}l^{-\mu+\nu+2}n_0^{\mu-\nu}$ a new parameter, and $Q_c$ a critical value depending on the exponents $\mu$ and $\nu$ and the initial distribution parameter $\alpha$. More importantly, a condition for $Q$ which holds generally for our model in case of conditional persistence ($\mu \geq \nu$) is derived. In Section \ref{sec:num}, we propose a finite-difference algorithm to determine numerical solutions of the considered model. In Section \ref{sec:results}, we will make use of the  numerical algorithm presented in Section \ref{sec:num} to examine the influence of the exponents $\mu$ and $\nu$ as well as the initial distribution parameter on the critical value $Q_c(\mu,\nu,\alpha)$. This algorithm will be used to determine the minimum value of the initial distribution parameter that guarantees the persistence of the population. Moreover, we examine the conditions for population survival on either the habitat size (for a given total population) or the total population (for a given habitat size). Furthermore, a detail discussion of the outlined procedures and obtained results, accompanied by illustrative examples, is provided.  Finally, in Section \ref{sec:con} we present concluding remarks.

\section{The model\label{sec:model}}

The reference framework for our analysis is the reaction-diffusion equation (\ref{eq1}) introduced in \cite{colombo2018nonlinear}. 
The main properties of (\ref{eq1}) as a dynamical population model, along with its ecological implications, are as follows:
\begin{itemize}
\item The diffusion term $(D\,u^{\nu-1}\,u_x)_x$ is, in general, nonlinear. Nonlinear diffusion terms, particularly power-law diffusion coefficients, have been widely used in dynamical population models. On the one hand, the case $\nu>1$ corresponds to species that, to avoid extinction, exhibit a low dispersal rate at low concentrations or 
to species provided with regulatory mechanisms to reduce overcrowding: large populations tend to a faster dispersion~\cite{dornelas2019single,birzu2019genetic,newman1980some}. On the other hand, social species such as insects tend to cluster as a survival strategy~\cite{sumpter2006principles,deneubourg1991dynamics}, which corresponds to $0<\nu<1$. This behavior is also observed in invading species~\cite{shigesada1997biological}.
Incidentally, nonlinear diffusion processes are also relevant in other contexts, such as transport in porous media~\cite{vazquez2007porous}, interstellar diffusion~\cite{newman1981galactic} or energy spreading in nonlinear disordered lattices~\cite{mulansky2013energy}.
\item The reaction term $a\,u^{\mu}$ corresponds to a growth rate depending on the population density $a\,u^{\mu-1}$, which is particularly relevant at low concentrations. If $\mu>1$, the term describes a population that, at low densities, experiences difficulties in reproduction and grows very slowly. This phenomenon is known as Allee effect~\cite{courchamp2008allee}. The opposite behaviour, $\mu<1$, is also possible in species, for instance zooplankton, which are capable to fit their reproduction mechanisms depending on the population density~\cite{gerritsen1977encounter}.
\item The intraspecific competition term has been neglected. Typical reasons to remove this term~\cite{murray2007mathematical} include: presence of abundant resources, human intervention (for instance, population in an experimental environment where nutrients are being continuously supplied) or social animal populations, where the cooperative behavior may significantly reduce the effect of the intraspecific competition for resources.
\end{itemize}
Equation~(\ref{eq1}) with hostile (Dirichlet) boundary conditions~(\ref{boun}) and the initial condition $u(x,0)\,=\,u_0 (x)$ corresponding to a total population $n_0$~(\ref{n0}) leads us to the problem~(\ref{prob1})-(\ref{n0}). Let us start by rewriting problem~(\ref{prob1})-(\ref{n0}) in terms of nondimensional variables. If we rescale $x$ by using the habitat size, i.e.
\begin{equation}\label{nvar11}
X\,=\,\frac{x}{l},
\end{equation}
it is possible to remove $a$ and $D$ from the reaction-diffusion equation provided that $\mu\,\neq\,\nu$, by rescaling $t$ and $u$ as
\begin{equation}\label{nvar12}
T\,=\,a\,\left(\frac{D}{a\,l^2}\right)^{\frac{\mu-1}{\mu-\nu}}\,t,\quad
\rho\,=\,\left(\frac{D}{a\,l^2}\right)^{-\frac{1}{\mu-\nu}}\,u.
\end{equation}
In this case, problem~(\ref{prob1}) takes the form
\begin{equation}\label{prob2}
\left\{\begin{array}{ll}
\rho_T\,=\,(\rho^{\nu-1}\,\rho_X)_X\,+\,\rho^{\mu}, \quad & X\,\in\,\left[-\frac{1}{2},\frac{1}{2}\right],\; T\,>\,0,\\  \\
\rho\left(-\frac{1}{2},T\right)\,=\,\rho\left(\frac{1}{2},T\right)\,=\,0, & T\,>\,0,\\ \\
\rho(X,0)\,=\,\rho_0(X), & X\,\in\,\left[-\frac{1}{2},\frac{1}{2}\right],
\end{array}\right.\end{equation}
where
\begin{equation}\label{rho01}
\rho_0(X)\,=\,\left(\frac{D}{a\,l^2}\right)^{-\frac{1}{\mu-\nu}}\,u_0\left(l\,X\right),\end{equation}
which satisfies
\begin{equation}\label{rho0N0}
\int_{-\frac{1}{2}}^{\frac{1}{2}}\rho_0(X)\,dX\,=\,N_0,
\end{equation}
with
\begin{equation}\label{N0n0}
N_0\,=\,\left(\frac{D}{a\,l^2}\right)^{-\frac{1}{\mu-\nu}}\,\frac{n_0}{l}.
\end{equation}
The survival condition for~(\ref{prob2}) with $\rho_0(X)$ satisfying~(\ref{rho0N0}) takes the form
$N_0\,\geq\,N_{0,c}$
where $N_{0,c}$ is a critical value depending on the growth exponent $\mu$, the diffusion exponent $\nu$ and the particular choice of the initial condition $\rho_0(X)$ satisfying~(\ref{rho0N0}). If we assume that we are considering a family of initial data satisfying~(\ref{rho0N0})
and depending on a certain set of parameters $\Lambda\,=\,\{\alpha_1,\alpha_2,\dots,\alpha_K\}$ the survival condition can be written as
\begin{equation}\label{cond1}
N_0\,\geq\,N_{0,c}(\mu,\nu,\Lambda).
\end{equation}
On the other hand, if $\nu\,=\,\mu$, in order to remove the parameters from the diffusion-reaction equation we take into account~(\ref{n0}) and use the variables
\begin{equation}\label{nvar21}
X\,=\,\sqrt{\frac{a}{D}}\,x,\quad 
T\,=\,a\,\left(\frac{a}{D}\right)^{\frac{\mu-1}{2}}\,n_0^{\mu-1}\,t, \quad
\rho\,=\,\left(\frac{a}{D}\right)^{-\frac{1}{2}}\frac{u}{n_0}.
\end{equation}
In terms of~(\ref{nvar21}), problem~(\ref{prob1}) takes the form
\begin{equation}\label{prob3}
\left\{\begin{array}{ll}
\rho_T\,=\,(\rho^{\mu-1}\,\rho_X)_X\,+\,\rho^{\mu}, \quad & X\,\in\,\left[-\frac{L}{2},\frac{L}{2}\right],\;  T\,>\,0,\\  \\
\rho\left(-\frac{L}{2},T\right)\,=\,\rho\left(-\frac{L}{2},T\right)\,=\,0, & T\,>\,0,\\  \\
\rho(X,0)\,=\,\rho_0(X), & X\,\in\,\left[-\frac{L}{2},\frac{L}{2}\right],
\end{array}\right.\end{equation}
where
\begin{equation}\label{L}
L\,=\,\sqrt{\frac{a}{D}}\,l,
\end{equation}
and
\begin{equation}\label{rho02}
\rho_0(X)\,=\,\sqrt{\frac{D}{a}}\,\frac{1}{n_0}\,u_0\left(\sqrt{\frac{D}{a}}X\right),
\end{equation}
that satisfies
\begin{equation}\label{N0eq2}
\int_{-\frac{L}{2}}^{\frac{L}{2}}\rho_0(X)\,dX\,=\,
1.
\end{equation}
Colombo and Anteneodo in~\cite{colombo2018nonlinear} provided numerical evidence that the survival condition for~(\ref{prob3}) is given by
\begin{equation}\label{cond2} 
L\,\geq\,\frac{\pi}{\sqrt{\mu}}, 
\end{equation}
where $L_c\,=\,\frac{\pi}{\sqrt{\mu}}$ is the value of $L$ corresponding at the unique nonnegative steady solution of~(\ref{prob3}), i.e.
$$\rho_s(X)\,=\,C\,\left(\cos(\sqrt{\mu}\,X)\right)^{\frac{1}{\mu}},$$
with $C$ being a positive constant depending on the initial data $\rho_0(X)$.
Next, we note that both survival conditions~(\ref{cond1}) and~(\ref{cond2}) can be rewritten in a unique form. Indeed, if~(\ref{N0n0}) is replaced into~(\ref{cond1}) and the two sides of the resulting inequality are raised to $\mu-\nu$ (recall that $\mu\,>\,\nu$), and~(\ref{L}) is replaced into~(\ref{cond2}) and both sides of the resulting inequality are raised to 2, then,  conditions~(\ref{cond1}) and~(\ref{cond2}) are reformulated as
\begin{equation}\label{fcond}  Q\,\geq\,Q_c(\mu,\nu,\Lambda), \end{equation}
where $Q$ is the parameter defined in~(\ref{Q}) and 
\begin{equation}\label{Qc}\everymath{\displaystyle}
Q_c(\mu,\nu,\Lambda)\,=\,\left\{\begin{array}{ccc}
N_{0,c}(\mu,\nu,\Lambda)^{\mu-\nu} & \mbox{if} & \mu\,>\,\nu,\\  \\
\frac{\pi^2}{\mu}                  & \mbox{if} & \mu\,=\,\nu.
\end{array}\right.
\end{equation}
In terms of $Q$, the equation~(\ref{rho0N0}) for the case $\mu\,>\,\nu$ reads
\begin{equation}\label{rho0Q}
\int_{-\frac{1}{2}}^{\frac{1}{2}}\rho_0(X)\,dX\,=\,Q^{\frac{1}{\mu-\nu}}.
\end{equation}
Some consequences of~(\ref{Q}), (\ref{fcond}) are:
\begin{enumerate}
\item For the particular case $\mu\,=\,\nu\,+\,2$, the equations~(\ref{Q}), (\ref{fcond}) do not impose any conditions on the habitat size $l$, but only 
on the total initial population $n_0$. Indeed, if $\mu\,=\,\nu\,+\,2$ (and consequently $\mu\,>\,2$)~(\ref{Q}), (\ref{fcond}) reads
\begin{equation}\label{critcase}
n_0\,\geq\,\sqrt{\frac{D}{a}Q_c(\mu,\mu-2,\Lambda)} \, .
\end{equation}
\item For a given total population $n_0$ and $\mu\,\in\,(\nu,\,\nu\,+\,2)$ the equations~(\ref{Q}), (\ref{fcond}) impose a condition
on the habitat size of the form
\begin{equation}\label{min}
l\,\geq\,l_c,
\end{equation}
where
\begin{equation}\label{critsize}
l_c\,=\,\left(\frac{D}{a}Q_c(\mu,\nu,\Lambda)\right)^{\frac{1}{-\mu+\nu+2}}n_0^{\frac{\mu-\nu}{\mu-\nu-2}},
\end{equation}
i.e., the population survives supposed that the habitat size is bigger than the critical size~(\ref{critsize}).
\item For a given total population $n_0$ and $\mu\,>\,\nu\,+\,2$ (and consequently $\mu\,>\,2$), the equations~(\ref{Q}), (\ref{fcond}) impose a condition on the habitat size of the form
\begin{equation}\label{max}
l\,\leq\,l_c,
\end{equation}
with $l_c$ being given by~(\ref{critsize}), i.e., the population survives supposed that the habitat size is smaller than the critical size~(\ref{critsize}).
\end{enumerate}
Note that these conclusions are in agreement with the results in~\cite{colombo2018nonlinear}. However, we extend the analysis by applying a numerical algorithm to investigate for the families of initial data~(\ref{u01}) and~(\ref{u02}), the dependence of the critical values, which we conveniently denote by $Q_{c,1}(\mu,\nu,\alpha)$ and $Q_{c,2}(\mu,\nu,\alpha)$, respectively, not only with respect to the exponents $\mu$ and $\nu$, but also with respect to the parameter $\alpha$ which governs the shape of the initial data. Furthermore, we employ the numerical algorithm to determine the required value of $\alpha$ to ensure the survival of the population.

\section{The numerical algorithm\label{sec:num}}

In order to present the numerical algorithm to approximate the solutions of problems~(\ref{prob2}) and~(\ref{prob3}) we consider the more general problem 
\begin{equation}\label{probg}
\left\{\begin{array}{ll }
\rho_T\,=\,(\rho^{\nu-1}\,\rho_X)_X\,+\,\rho^{\mu}, \quad &  X\,\in\,\left[-\frac{L}{2},\frac{L}{2}\right],\; T\,>\,0,\\  \\
\rho\left(-\frac{L}{2},T\right)\,=\,\rho\left(-\frac{L}{2},T\right)\,=\,0, &  T\,>\,0,\\  \\
\rho(X,0)\,=\,\rho_0(X), &  X\,\in\,\left[-\frac{L}{2},\frac{L}{2}\right].
\end{array}\right.\end{equation}
Note that problem~(\ref{prob2}) is a particular case of~(\ref{probg}) for $L=1$ while~(\ref{prob3}) is a particular case of~(\ref{probg}) with $\nu\,=\,\mu$. Let the solution domain of~(\ref{probg}) be the closed subset of $\mathbb{R}^2$ given by
$$ \Omega= \left[ -\frac{L}{2},\frac{L}{2} \right] \times  \left[ 0,\tilde{T} \right],$$
with $\tilde{T}>0$. Let subdivide $\Omega$ into uniform rectangular meshes formed by the intersection points, usually called grid points, of the lines $X_i=-\frac{L}{2}+ i h$, $i=0(1)m$, $T_j=j k$, $j=0(1)n$, where $m h=L$, $n k=\tilde{T}$.

PDE in (\ref{probg}) can be written equivalently as follows
\begin{equation}
\label{prob4}
\rho_T= \frac{1}{\nu}\left( \rho^{\nu} \right)_{XX}+ \rho^{\mu}.
\end{equation}
We will replace the derivatives of PDE (\ref{prob4}) on each grid point by a finite-difference approximation in terms of the values of $\rho$ at neighbouring grid points. Here, we consider the finite-difference representation of PDE (\ref{prob4}) given by
\begin{equation}
\label{numscheme1} \frac{\rho_{i,j+1}-\rho_{i,j}}{k}=\frac{1}{2 \nu} \Big( \delta_X^2  \rho_{i,j+1}^\nu + \delta_X^2 \rho_{i,j}^\nu  \Big)+\rho_{i,j}^\mu,
\end{equation}
which has been obtained by employing a forward-difference for $\rho_T$, the Crank-Nicolson approximation to $\left( \rho^{\nu} \right)_{XX}$, and then making use of the usual central-difference operator $\delta_X^2$ defined by
$$ \delta_X^2 \rho_{i,j}= \frac{\rho_{i+1,j}-2\rho_{i,j}+\rho_{i-1,j}}{h^2}.$$

Clearly, the finite-difference scheme (\ref{numscheme1}) is nonlinear. Therefore, in order to linearize scheme (\ref{numscheme1}) we apply the Richtmyer's linearization method \cite{richtmyer1967,smith1985}. Let us consider the Taylor's expansion of $\rho_{i,j+1}^\nu$ about the point $(i,j)$
$$\begin{array}{rcl}
\rho_{i,j+1}^\nu &=& \rho_{i,j}^\nu+ k \dfrac{\partial \rho_{i,j}^\nu}{\partial t}+\cdots \vspace*{0.2cm}\\
&=& \rho_{i,j}^\nu+ k \nu \rho_{i,j}^{\nu-1}\dfrac{\partial \rho_{i,j}}{\partial t}+\cdots
\end{array}$$
Hence to terms of order $k$ 
\begin{equation}\label{approx1}
\begin{array}{rcl}
\rho_{i,j+1}^\nu &\cong& \rho_{i,j}^\nu+ k \nu \rho_{i,j}^{\nu-1} \left(\dfrac{\rho_{i,j+1}-\rho_{i,j}}{k}\right) \vspace*{0.2cm}\\
&=& \rho_{i,j}^\nu+  \nu \rho_{i,j}^{\nu-1} \left(\rho_{i,j+1}-\rho_{i,j}\right),
\end{array}
\end{equation}
we can substitute the nonlinear unknown $\rho_{i,j+1}^\nu$ by a linear approximation in $\rho_{i,j+1}$. Similarly,
\begin{equation}\label{approx2}
\begin{array}{rcl}
\rho_{i+1,j+1}^\nu &\cong& \rho_{i+1,j}^\nu+ \nu \rho_{i+1,j}^{\nu-1} \left(\rho_{i+1,j+1}-\rho_{i+1,j}\right) \vspace*{0.2cm}\\
\rho_{i-1,j+1}^\nu &\cong& \rho_{i-1,j}^\nu+ \nu \rho_{i-1,j}^{\nu-1} \left(\rho_{i-1,j+1}-\rho_{i-1,j}\right)
\end{array}
\end{equation}
Substituting (\ref{approx1}) and (\ref{approx2}) into scheme (\ref{numscheme1}) and defining $W_i=\rho_{i,j+1}-\rho_{i,j}$, after some simplifications, we obtain a linear system for $W_i$
\begin{equation}
\label{numscheme2} \begin{array}{l} \displaystyle \rho_{i+1,j}^{\nu-1}W_{i+1}-2\left(\rho_{i,j}^{\nu-1}+\frac{h^2}{k}\right)W_{i}+\rho_{i-1,j}^{\nu-1}W_{i-1} \vspace*{0.2cm}\\
\displaystyle \qquad \qquad \quad =-\frac{2}{\nu}\rho_{i+1,j}^{\nu} + \frac{4}{\nu}\rho_{i,j}^{\nu} -\frac{2}{\nu}\rho_{i-1,j}^{\nu}- 2h^2\rho_{i,j}^\mu,
\end{array}\end{equation}
$i=1(1)m-1$. The solution at the $(j+1)$th time-level is obtained from $\rho_{i,j+1}=W_i+\rho_{i,j}$.

From boundary conditions given in (\ref{probg}), we have $\rho_{0,j}=\rho_{m,j}=0$, $\forall j$, and therefore $W_0=W_m=0$. Let 
$$\begin{array}{rcl}
\bm{{\rho^{k}}}^{(j)}&=& \left(\rho_{1,j}^{k},\rho_{2,j}^{k}, \ldots, \rho_{m-1,j}^{k}\right)^T,\vspace*{0.2cm}\\
\bm{W}&=& \left(W_{1},W_{2}, \ldots, W_{m-1}\right)^T.
\end{array}
$$
Consequently, equation (\ref{numscheme2}) can be written in matrix form
\begin{equation}\label{numscheme3} A \bm{W}= -\frac{2}{\nu}B \bm{{\rho^{\nu}}}^{(j)}-2 h^2 \bm{{\rho^{\mu}}}^{(j)},\end{equation}
where $A=B \,\mbox{diag}\{\bm{{\rho^{\nu-1}}}^{(j)}\}- \frac{2 h^2}{k} I_{m-1}$ is a tridiagonal matrix, $I_{m-1}$ is the unit matrix of order $m-1$ and $B$ is given by 
$$B=\left(%
\begin{array}{ccccc}
  -2 & 1 &  &  &     \\
  1 & -2 &  1 &  &     \\
   &  \ddots & \ddots & \ddots    \\
   &  &   1 & -2 &  1  \\
   &  &    & 1 & -2  \\
\end{array}%
\right).$$
Matrix system (\ref{numscheme3}) can be solved in an efficient way by using Thomas algorithm \cite{contedeboor1972}.

\section{Results and discussion \label{sec:results}}
In this section, we will make use of the numerical algorithm presented in Section \ref{sec:num} to determine the dependence of the critical values $Q_{c,1}(\mu,\nu,\alpha)$ and $Q_{c,2}(\mu,\nu,\alpha)$ on the values of the exponents $\mu$, $\nu$ in the model and the parameter $\alpha$ in the initial conditions~(\ref{u01}) and~(\ref{u02}), respectively. Furthermore, we will estimate the minimum value of $\alpha$ for which the survival of the population is guaranteed.

In terms of the nondimensional variables, the initial condition in~(\ref{prob2}) for the case $\mu\, > \, \nu$ is determined,  respectively, by~(\ref{u01}), (\ref{rho01}) and~(\ref{Q}) as
\begin{equation}\label{rho0fam1}\everymath{\displaystyle}
\rho_{0,1}(X;\alpha)\,=\, \frac{Q^{\frac{1}{\mu-\nu}}}{B\left(1+\alpha,1+\alpha\right)}\left(\frac{1}{4}\,-\,X^2\right)^{\alpha},
\end{equation}
or by~(\ref{u02}), (\ref{rho01}) and~(\ref{Q}) as
\begin{equation}\label{rho0fam2}\everymath{\displaystyle}
\rho_{0,2}(X;\alpha)\,=\, \frac{Q^{\frac{1}{\mu-\nu}}}{B(1+ \gamma(\alpha),1+2  \gamma(\alpha))}\left(\left(\frac{1}{2}-X\right)\left(\frac{1}{2}+X\right)^2\right)^{ \gamma(\alpha)},
\end{equation}
both satisfying~(\ref{rho0Q}).

Analogously, the initial conditions in~(\ref{prob3}) for the case $\mu \, = \, \nu$ are given, respectively, by
\begin{equation}\label{rho0f1eq}
\rho_{0,1}(X;\alpha)\,=\,\frac{Q^{-\frac{1}{2}}}{B(1+\alpha,1+\alpha)}\left(\frac{1}{4}-\frac{X^2}{Q}\right)^{\alpha},
\end{equation}
and
\begin{equation}\label{rho0f2eq}
\rho_{0,2}(X;\alpha)\,=\,\frac{Q^{-\frac{1}{2}}}{B(1+\gamma(\alpha),1+2\gamma(\alpha))}\left(\left(\frac{1}{2}-\frac{X}{Q^{\frac{1}{2}}}\right)\left(\frac{1}{2}+\frac{X}{Q^{\frac{1}{2}}}\right)^2\right)^{\gamma(\alpha)},
\end{equation}
both satisfying~(\ref{N0eq2}).

\subsection{Population density and total population dynamics}

Let us start by showing some illustrative examples of the behaviour of the population density and the total population dynamics described by the solutions of~(\ref{prob2}) in relation to the critical values $Q_{c,1}(\mu,\nu,\alpha)$ and $Q_{c,2}(\mu,\nu,\alpha)$ for the families of initial conditions (\ref{rho0fam1}) and (\ref{rho0fam2}), respectively.

We estimate the total population (in the nondimensional variables)
\begin{equation}\label{NT}
N(T)\,=\,\int_{- \frac{1}{2}}^{\frac{1}{2}}\rho(X,T)\,dX,
\end{equation}
by using the trapezoidal rule as
\begin{equation}\label{NTapprox}
N(j\,k)\,\approx\,h\,\left(\frac{\rho_{0,j}}{2}\,+\,\sum_{i=1}^{n-1}\rho_{i,j}\,+\,\frac{\rho_{n,j}}{2}\right)\quad j\,=\,0,\,1,\,\dots,\,m.
\end{equation}

\begin{figure}[tbp]
	\centering
        \includegraphics[width=12cm]{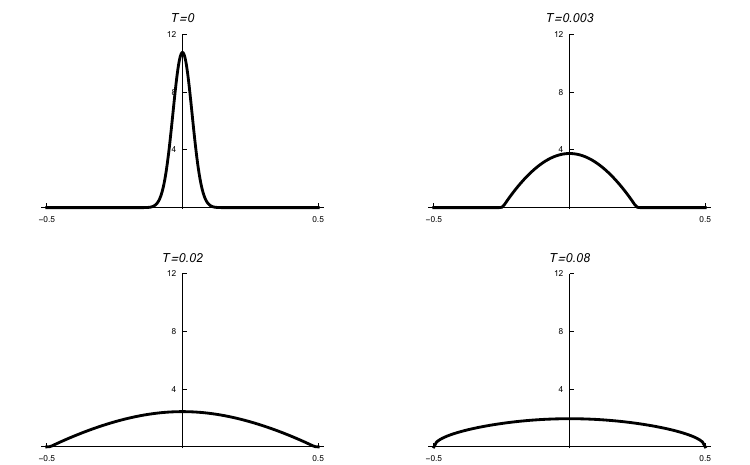}
	\caption{Evolution of the population density governed by (\ref{prob2}) in the nonlinear case $\mu=4$, $\nu=2$, when $Q=0.9$ and initial condition (\ref{rho0fam1}) with $\alpha=100$ is considered. 
	\label{fig:tempevol1ext}}
\end{figure}

\begin{figure}[tbp]
	\centering
        \includegraphics[width=12cm]{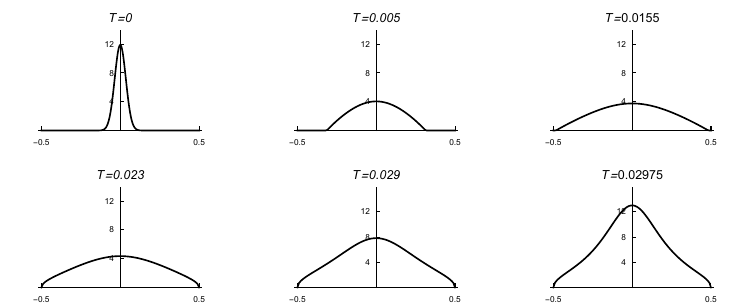}
	\caption{Evolution of the population density governed by (\ref{prob2}) in the nonlinear case $\mu=4$, $\nu=2$, and when $Q=1.1$ and initial condition (\ref{rho0fam1}) with $\alpha=100$ is considered. 
	\label{fig:tempevol1crec}}
\end{figure}

\begin{figure}[tbp]
	\centering
        \includegraphics[width=12cm]{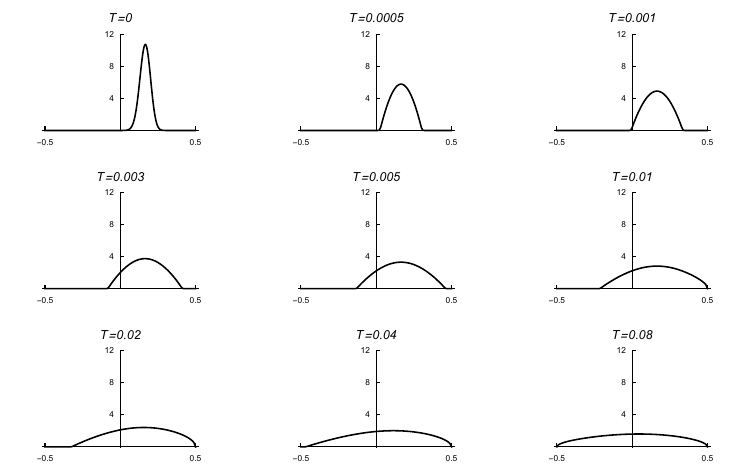}
	\caption{Evolution of the population density governed by (\ref{prob2}) in the nonlinear case $\mu=4$, $\nu=2$, when $Q=0.9$ and initial condition (\ref{rho0fam2}) with $\alpha=100$ is considered. 
	\label{fig:tempevol2ext}}
\end{figure}

\begin{figure}[tbp]
	\centering
        \includegraphics[width=12cm]{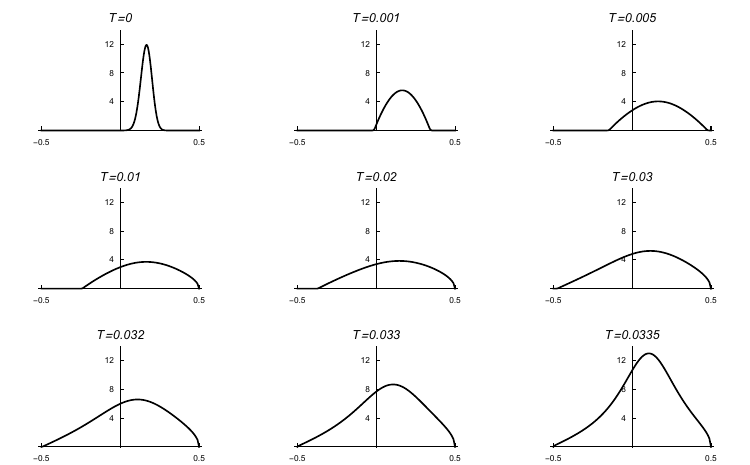}
	\caption{Evolution of the population density governed by (\ref{prob2}) in the nonlinear case $\mu=4$, $\nu=2$, when $Q=1.1$ and initial condition (\ref{rho0fam2}) with $\alpha=100$ is considered. 
	\label{fig:tempevol2crec}}
\end{figure}

From Section \ref{sec:model} we know that for a value of $Q$ less than the corresponding critical value, the population becomes extinct (see Figures \ref{fig:tempevol1ext} and \ref{fig:tempevol2ext}), while for $Q$ greater than the corresponding critical value, increases without limit (see Figures \ref{fig:tempevol1crec} and \ref{fig:tempevol2crec}). In order to illustrate this situation, we focus on the nonlinear case $\mu=4$, $\nu=2$ and $\alpha=100$. 

We observe that when the initial distribution is symmetric (\ref{rho0fam1}), the population density rapidly declines symmetrically, keeping its maximum value in the midpoint of the habitat, $X=0$. After an interval of time, the situation drastically changes depending on whether $Q<Q_{c,1}(4,2,100)$ or $Q>Q_{c,1}(4,2,100)$. If $Q<Q_{c,1}(4,2,100)$, the population density tends to $0$ at all points of the habitat size, which leads to the extinction of the given population (see Figure \ref{fig:tempevol1ext}). If $Q>Q_{c,1}(4,2,100)$, the population density start growing symmetrically with respect to the midpoint of the habitat without bound (see Figure \ref{fig:tempevol1crec}).

Analogously, we observe that when the initial distribution is asymmetric (\ref{rho0fam2}), the population density at the original maximum location, $X=1/6$, rapidly declines since individuals spread out across the entire habitat size in search of less crowded conditions. Once the individuals have occupied nearly the entire habitat size, the situation drastically changes depending on whether $Q<Q_{c,2}(4,2,100)$ or $Q>Q_{c,2}(4,2,100)$. If $Q<Q_{c,2}(4,2,100)$, the population density tends to zero exhibiting a behavior similar to that observed when initial distribution (\ref{rho0fam1}) is considered with $Q<Q_{c,1}(4,2,100)$ (see Figure \ref{fig:tempevol2ext}). If $Q>Q_{c,2}(4,2,100)$, the population density grows without bound (see Figure \ref{fig:tempevol2crec}).

\begin{figure}[tbp]
	\centering
        \includegraphics[width=10cm]{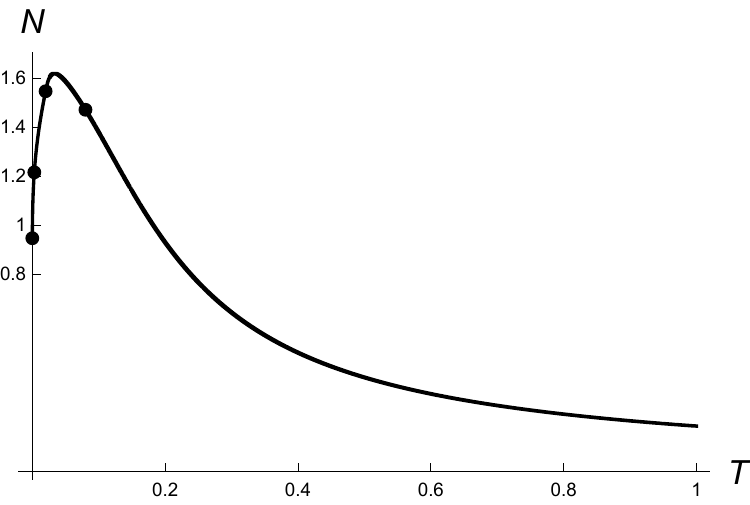}
	\caption{Evolution of the total population governed by~(\ref{prob2}) in the nonlinear case $\mu=4$, $\nu=2$, when $Q=0.9$ and initial condition (\ref{rho0fam1}) with $\alpha=100$ is considered. Black points represent the total population size for the temporal values considered in Figure \ref{fig:tempevol1ext}.
	\label{fig:Nttempevol1ext}}
\end{figure}

\begin{figure}[tbp]
	\centering
        \includegraphics[width=10cm]{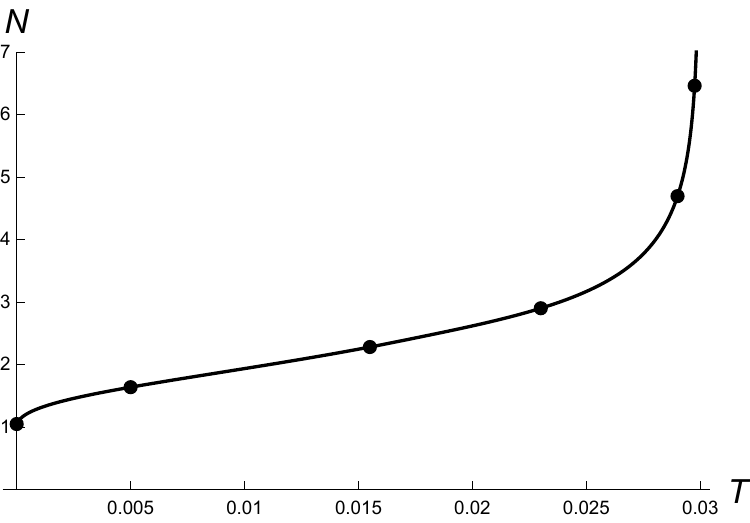}
	\caption{Evolution of the total population governed by~(\ref{prob2}) in the nonlinear case $\mu=4$, $\nu=2$, when $Q=1.1$ and initial condition (\ref{rho0fam1}) with $\alpha=100$ is considered. Black points represent the total population size for the temporal values considered in Figure \ref{fig:tempevol1crec}.
	\label{fig:Nttempevol1crec}}
\end{figure}

\begin{figure}[tbp]
	\centering
        \includegraphics[width=10cm]{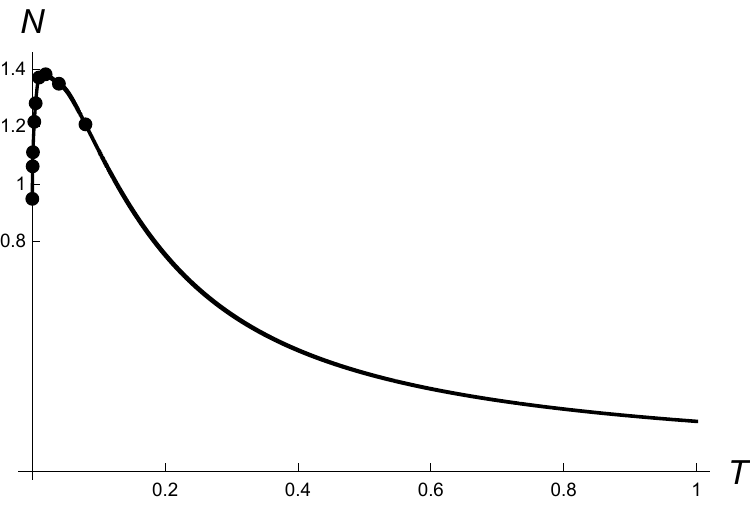}
	\caption{Evolution of the total population governed by~(\ref{prob2}) in the nonlinear case $\mu=4$, $\nu=2$, when $Q=0.9$ and initial condition (\ref{rho0fam2}) with $\alpha=100$ is considered. Black points represent the total population size for the temporal values considered in Figure \ref{fig:tempevol2ext}.
	\label{fig:Nttempevol2ext}}
\end{figure}

\begin{figure}[tbp]
	\centering
        \includegraphics[width=10cm]{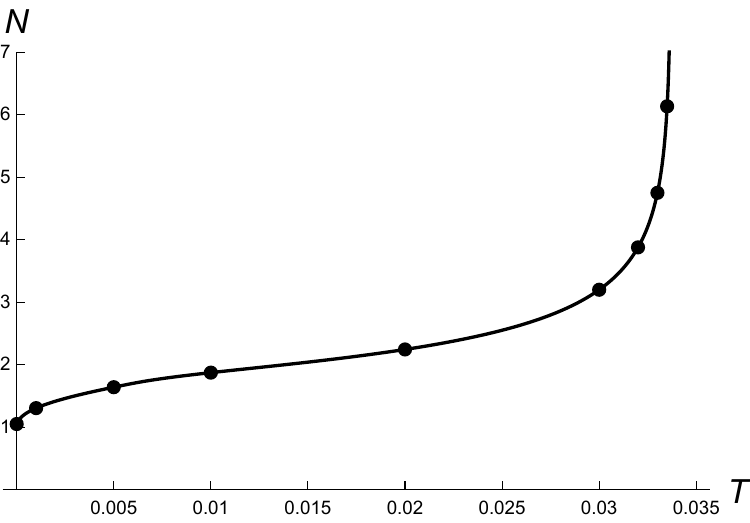}
	\caption{Evolution of the total population governed by~(\ref{prob2}) in the nonlinear case $\mu=4$, $\nu=2$, when $Q=1.1$ and initial condition (\ref{rho0fam2}) with $\alpha=100$ is considered. Black points represent the total population size for the temporal values considered in Figure \ref{fig:tempevol2crec}.
	\label{fig:Nttempevol2crec}}
\end{figure}

Now, in Figures \ref{fig:Nttempevol1ext}, \ref{fig:Nttempevol1crec}, \ref{fig:Nttempevol2ext} and \ref{fig:Nttempevol2crec}, we show the temporal evolution of the total population governed by~(\ref{prob2}) for the different cases previously considered in Figures \ref{fig:tempevol1ext}, \ref{fig:tempevol1crec}, \ref{fig:tempevol2ext} and \ref{fig:tempevol2crec}, where black points represent the total population size for the temporal values considered in the above mentioned figures.

\subsection{The critical values $Q_{c,1}(\mu,\nu,\alpha)$ and $Q_{c,2}(\mu,\nu,\alpha)$: numerical estimation}\label{sec:critval}

Our first goal is to apply the numerical algorithm to estimate $Q_{c,1}(\mu,\nu,\alpha)$ and $Q_{c,2}(\mu,\nu,\alpha)$ for different values of the exponents $\mu$, $\nu$ in the model and the parameter
$\alpha$, respectively, in the initial conditions~(\ref{rho0fam1}) and~(\ref{rho0fam2}).  

For simplicity of notation, we outline the procedure for estimating the critical values by considering the first family of initial conditions. Let us start by setting particular values of $\mu_0$, $\nu_0$, $\alpha_0$ and apply the numerical algorithm to~(\ref{prob2}) where we choose $Q=Q_{0}$, for a certain $Q_{0}$ for the family of initial conditions ~(\ref{rho0fam1}). 

We have seen that for small enough $Q_0$ the approximated total population is an asymptotically decreasing function describing the extinction of the population, while for large enough
$Q_0$ the total population after an interval of time turns to be a monotonically increasing function that reaches arbitrarily large values (with possible blow-up). Next, we let $Q$ take values until we obtain two closed enough values $Q^*$ and $Q^{**}$ such that the approximated total population~(\ref{NTapprox}) tends to zero for $Q=Q^*$ and grows without limit for $Q=Q^{**}$. To this end, we start with $Q=Q^{(0)}$ sufficiently large so that the total population becomes an increasing function once enough time has elapsed, choose a step size $\Delta Q$ sufficiently small compared with $Q^{(0)}$ and consider the sequence of values $Q^{(r)}=Q^{(r-1)}-\Delta Q$, for $r \in \mathbb{N}$. For increasing values of $r$, we apply the numerical algorithm to approximate the total population~(\ref{NTapprox}) corresponding to the initial condition~(\ref{rho0fam1}) with $Q=Q^{(r)}$. This process continues until, for a certain $R$, the approximated total population becomes an asymptotically decreasing function which describes the extinction of the population.  We denote $Q^*=Q^{(R)}$ and $Q^{**}=Q^{(R-1)}$. Then, we approximate the critical value
$Q_{c,1}(\mu_0,\nu_0,\alpha_0)\approx(Q^*+Q^{**})/2$. The same method is applied to the second family of initial conditions.

\begin{figure}[H]
	\centering
\includegraphics[width=6.6cm]{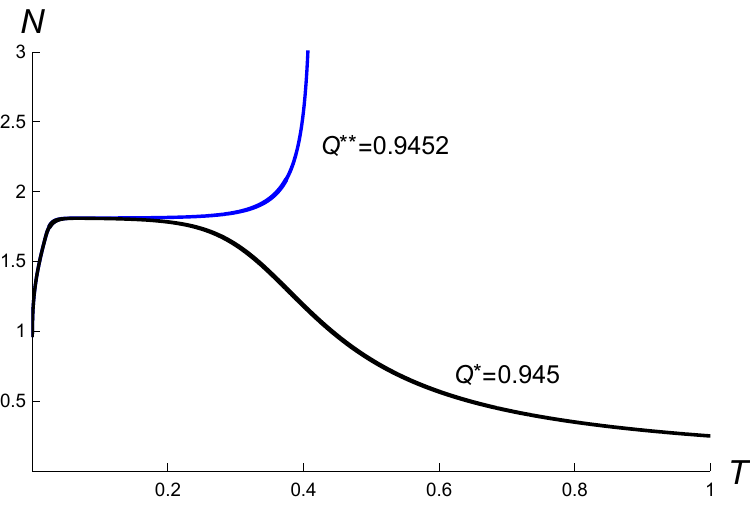}\hspace{2mm}
\includegraphics[width=6.6cm]{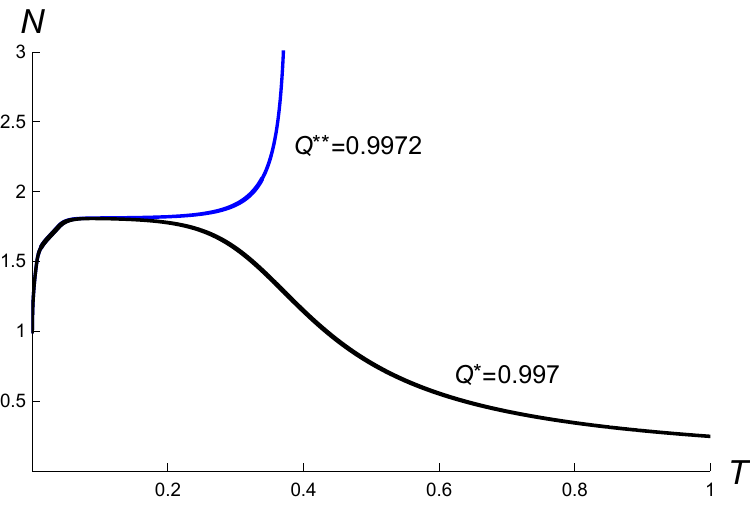}\\
	\caption{Evolution of the total population for equation~(\ref{prob2}) with $\mu_0=4$, $\nu_0=2$, $\alpha_0=100$ and initial condition~(\ref{rho0fam1}) (left) and~(\ref{rho0fam2}) (right). The left plot shows that the total population grows for $Q^{**}=0.9452$ (blue line) while it decreases for $Q^{*}=0.945$ (black line). We estimate the critical value $Q_{c,1}(4,2,100) \approx 0.9451$. Similarly, from the right plot it is observed how the total population grows for $Q^{**}=0.9972$ (blue line) while it decreases for $Q^*=0.997$ (black line). We estimate the critical value $Q_{c,2}(4,2,100) \approx 0.9971$.
	\label{fig:ejemQc}}
\end{figure}

We illustrate the method for $\mu_0=4$, $\nu_0=2$ and $\alpha_0=100$. The choice $\Delta Q=0.0002$ results in $Q^*=0.945$ and $Q^{**}=0.9452$ for the first family, and leads to $Q^*=0.997$ and $Q^{**}=0.9972$ for the second one. Thus, we approximate $Q_{c,1}(4,2,100) \approx 0.9451$ and $Q_{c,2}(4,2,100) \approx 0.9971$. Figure~\ref{fig:ejemQc} displays the approximate total population~(\ref{NTapprox}) for $Q^*$ (black lines) and $Q^{**}$ (blue lines). The left plot corresponds to the first family of initial conditions, whereas the right plot corresponds to the second family. 

\begin{figure}[H]
\centering
\includegraphics[width=6.6cm]{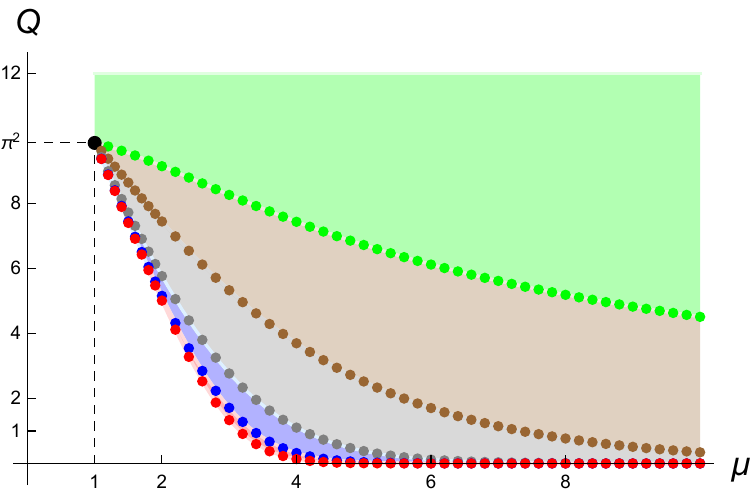}\hspace{2mm}
\includegraphics[width=6.6cm]{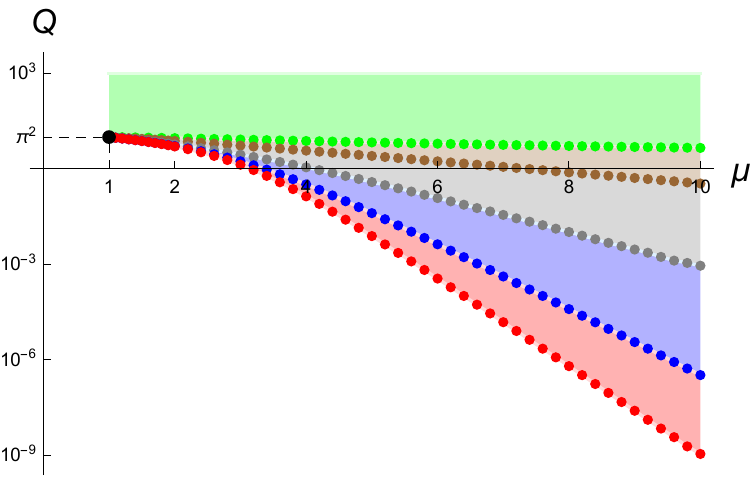}\\
\includegraphics[width=6.6cm]{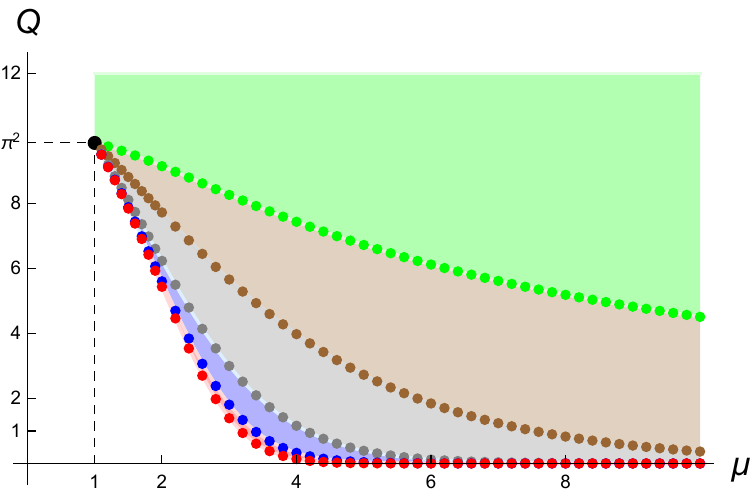}\hspace{2mm}
\includegraphics[width=6.6cm]{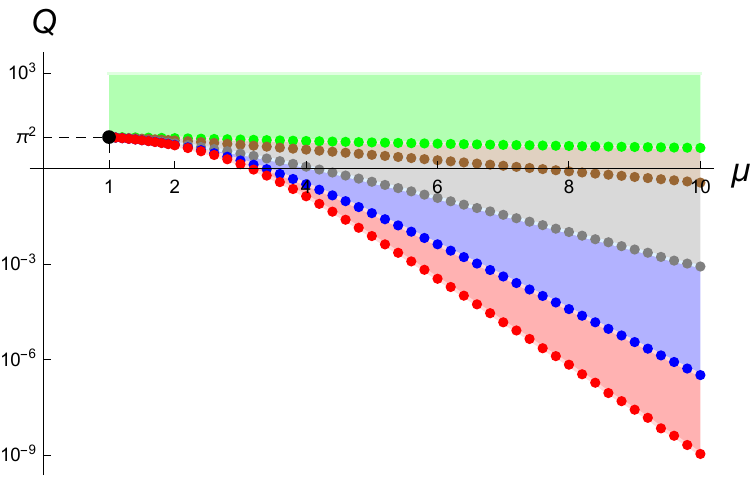}
\caption{Critical values $Q_{c,1}\left(\mu,1,\alpha_0\right)$ (top) and $Q_{c,2}\left(\mu,1,\alpha_0\right)$ (bottom) as functions of $\mu$ for $\alpha_0=0$ (green), 1 (brown), 10 (gray), 100 (blue) and 500 (red). The green shaded region corresponds to the survival region for $\alpha_0=0$ (i.e. the homogeneous initial distribution). As $\alpha_0$ increases, the survival region enlarges, adding progressively  the brown shaded region for $\alpha_0=1$, the gray shaded region for $\alpha_0=10$, the blue shaded region for $\alpha_0=100$ and the red shaded region for $\alpha_0=500$. The plots on the right are the same as those on the left, but using a logarithmic scale on the vertical axis to highlight the differences among the last three values of $\alpha_0$. The black point $(1,\pi^2)$ corresponds to $Q_{c,1}(1,1,\alpha)=Q_{c,2}(1,1,\alpha)=\pi^2$ for all $\alpha$.}
\label{fig:graphQnu1}
\end{figure}

Figure~\ref{fig:graphQnu1} displays the graphs of $Q_{c,1}(\mu,1,\alpha_0)$ (top) and $Q_{c,2}\left(\mu,1,\alpha_0\right)$ (bottom) as functions of $\mu$ for  $\alpha_0=0$ (green), $\alpha_0=1$ (brown), $\alpha_0=10$ (gray), $\alpha_0=100$ (blue) and $\alpha_0=500$ (red). Since $Q_{c,j}$, $j=1,2$, takes very small values as $\mu$ increases and $\alpha_0$ assumes larger values, a logarithmic scale is particularly appropriate to enhance visual resolution at smaller values while preserving proportional relationships whereas a linear scale is suitable for small values of $\mu$, where $Q$ spans a wide range of values. A linear scale is used in the plots on the left whereas a logarithmic scale is used in the plots on the right. The shaded regions correspond to the survival regions in the parameter plane $(\mu,Q)$. More precisely, the green shaded region is the survival region for 
the homogeneous initial distribution ($\alpha_0=0$) and as $\alpha_0$ increases the survival region enlarges progressively, adding the 
brown shaded region for $\alpha_0=1$, the gray shaded region for $\alpha_0=10$, the blue shaded region for $\alpha_0=100$ and the red shaded region for $\alpha_0=500$. We point out that when the scheme to estimate  $Q_{c,j}(\mu,1,\alpha_0)$ is applied to increasing values of $\mu$, the step size $\Delta Q$ may need to be reduced ensuring that it remains sufficiently small compared to the estimated value of 
$Q_{c,j}(\mu,1,\alpha_0)$. It can be appreciated that, for sufficiently large values of $\mu$  (approximately $\mu\gtrapprox4$), an increasing value of $\alpha_0$ significantly reduces the critical value $Q_{c,j}(\mu,1,\alpha_0)$ and consequently facilitates population survival. However, for smaller values of $\mu$ a considerably larger $\alpha$ is required to achieve a noticeable reduction in  $Q_c(\mu,1,\alpha_0)$, so that the strategy of using highly concentrated initial distributions becomes less effective. In particular, for $\mu=1$ (linear case), the critical value of $Q$ is independent of the initial distribution ($Q_{c,j}(1,1,\alpha)=\pi^2$, $j=1,2$).

\begin{figure}[H]
	\centering
        \includegraphics[width=13.5cm]{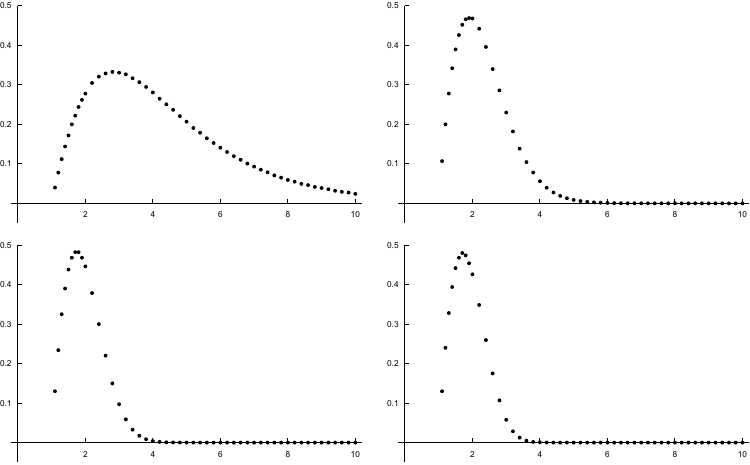}
        \caption{Differences $Q_{c,2}\left(\mu,1,\alpha_0\right)-Q_{c,1}\left(\mu,1,\alpha_0\right)$ between the critical values corresponding to the initial distributions~(\ref{rho0fam1}) and~(\ref{rho0fam2}) for $\alpha_0=1$ (top left), $\alpha_0=10$ (top right), $\alpha_0=100$ (bottom left) and $\alpha=500$ (bottom right).}\label{fig:graphdifQc}
        \end{figure}

Figure~\ref{fig:graphQnu1} also shows that the results are almost the same for both families of initial conditions. To provide a comparison between the critical values corresponding to initial conditions~(\ref{rho0fam1}) and~(\ref{rho0fam2}), we show in Figure~\ref{fig:graphdifQc} the differences $Q_{c,2}(\mu,1,\alpha_0)-Q_{c,1}(\mu,1,\alpha_0)$ for $\alpha_0=1$ (top left), $\alpha_0=10$ (top right), $\alpha_0=100$ (bottom left) and $\alpha_0=500$ (bottom right). Recall that $\alpha_0=0$ corresponds to the homogeneous distribution for both families and consequently the results are the same. 

Figure~\ref{fig:graphdifQc} illustrates that $Q_{c,2}(\mu,1,\alpha_0)\geq Q_{c,1}(\mu,1,\alpha_0)$ for all the cases we have analyzed. However, the differences are very small when compared with the values $Q_{c,2}(\mu,1,\alpha_0)$ or $Q_{c,1}(\mu,1,\alpha_0)$. That means that initial condition~(\ref{rho0fam1}) is slightly more propitious to ensure population survival. Note that the main difference between two families of initial conditions is that while~(\ref{rho0fam1}) are symmetric distributions centered at the origin, initial distributions~(\ref{rho0fam2}) are centered in $X=1/6$. Thus, under the initial distribution~(\ref{rho0fam1}), individuals take, on average, the same time to reach either endpoint of the interval. In contrast, individuals initially distributed according to~(\ref{rho0fam2}) tend to reach the right endpoint more quickly, which implies a shorter lifespan and a reduced period for reproduction. Conversely, those moving towards the left endpoint take longer to reach it, thereby enjoying a longer reproductive phase. Numerical results suggest that these opposing effects do not fully compensate each other, and the symmetric initial distribution~(\ref{rho0fam1}) appears to be slightly more favorable in terms of population survival.
 
\begin{figure}[H]
\centering
\includegraphics[width=6.6cm]{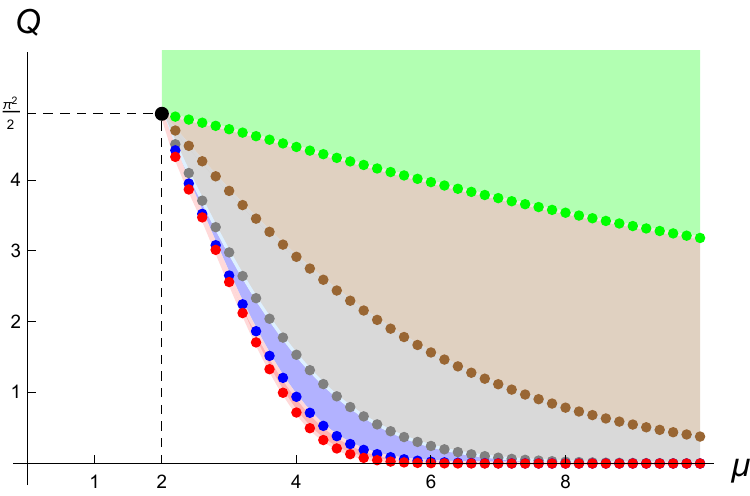}\hspace{2mm}
\includegraphics[width=6.6cm]{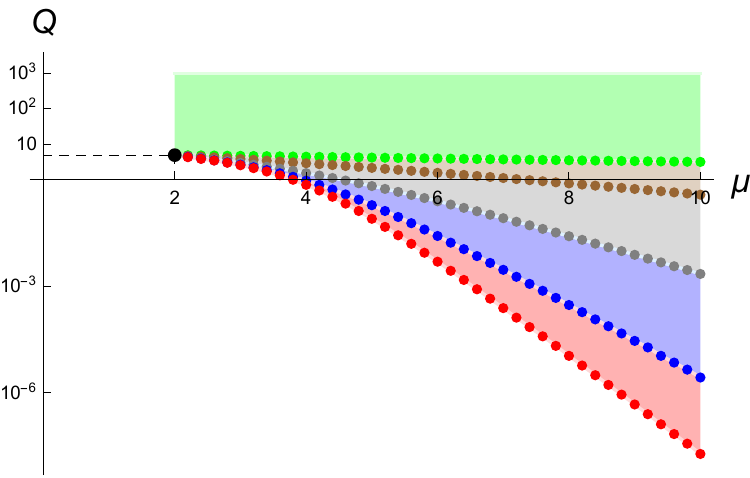}\\
\includegraphics[width=6.6cm]{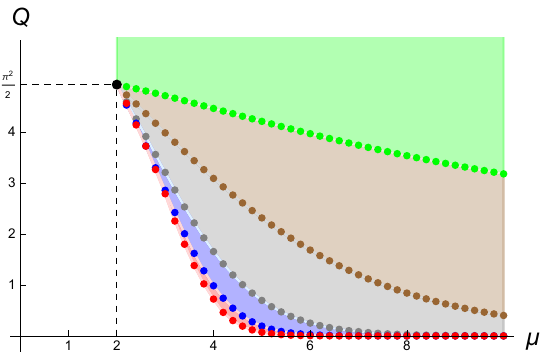}\hspace{2mm}
\includegraphics[width=6.6cm]{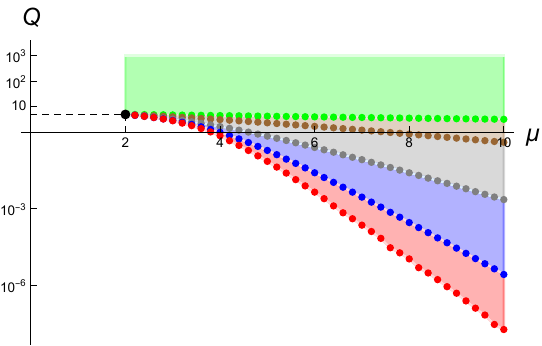}
\caption{Critical values $Q_{c,1}\left(\mu,2,\alpha_0\right)$ (top) and $Q_{c,2}\left(\mu,2,\alpha_0\right)$ (bottom) as functions of $\mu$ for $\alpha_0=0$ (green), 1 (brown), 10 (gray), 100 (blue) and 500 (red).  The green shaded region corresponds to the survival region for $\alpha_0=0$ (i.e. the homogeneous initial distribution). As $\alpha_0$ increases, the survival region enlarges, adding progressively  the brown shaded region for $\alpha_0=1$, the gray shaded region for $\alpha_0=10$, the blue shaded region for $\alpha_0=100$ and the red shaded region for $\alpha_0=500$.  The plots on the right are the same as those on the left, but using a logarithmic scale on the vertical axis to highlight the differences among the last three values of $\alpha_0$. The black point $(2,\pi^2/2)$ corresponds to $Q_{c,1}(2,2,\alpha)=Q_{c,2}(2,2,\alpha)=\pi^2/2$ for all $\alpha$.}
\label{fig:graphQnu2}
\end{figure}

Figure~\ref{fig:graphQnu2} shows the graphs of $Q_{c,1}(\mu,2,\alpha_0)$ (top) and $Q_{c,2}\left(\mu,2,\alpha_0\right)$ (bottom) as functions of $\mu$ using the same color code as in Figure~\ref{fig:graphQnu1}. Unlike the case $\nu_0=1$ illustrated in Figures~\ref{fig:graphQnu1} and \ref{fig:graphdifQc}, the case $\nu_0=2$ corresponds to a nonlinear diffusion term. Nevertheless, the
graphs can be interpreted in the same way as in Figure~\ref{fig:graphQnu1} for $\nu_0=1$.

\begin{figure}[H]
\centering
\includegraphics[width=6.6cm]{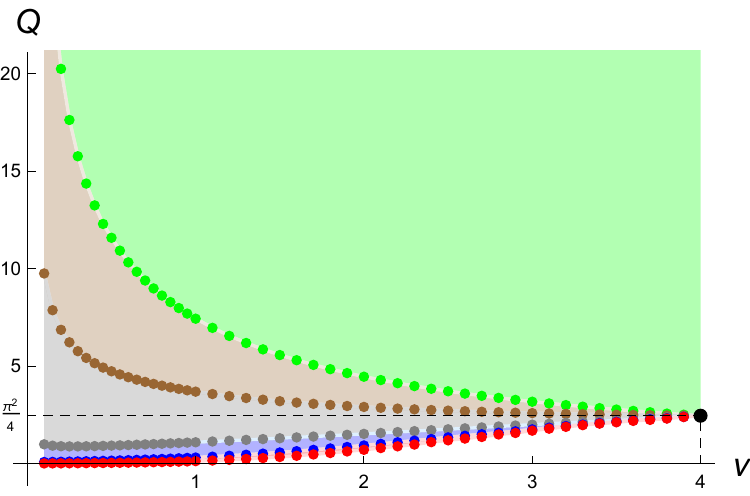}\hspace{2mm}
\includegraphics[width=6.6cm]{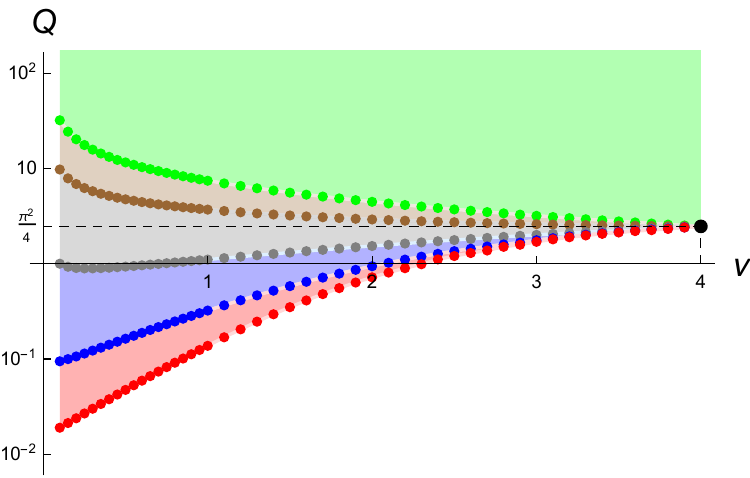}
\caption{Critical values $Q_{c,1}\left(4,\nu,\alpha_0\right)$ as functions of $\nu$ for $\alpha_0=0$ (green), 1 (brown), 10 (gray), 100 (blue) and 500 (red).  The green shaded region corresponds to the survival region for $\alpha_0=0$ (i.e. the homogeneous initial distribution). As $\alpha_0$ increases, the survival region enlarges, adding progressively  the brown shaded region for $\alpha_0=1$, the gray shaded region for $\alpha_0=10$, the blue shaded region for $\alpha_0=100$ and the red shaded region for $\alpha_0=500$.  The plot on the right is the same as the one on the left, but using a logarithmic scale on the vertical axis to highlight the differences among the last three values of $\alpha_0$. The black point $(4,\pi^2/4)$ corresponds to $Q_{c,1}(4,4,\alpha)=\pi^2/4$ for all $\alpha$.}
\label{fig:graphQmu4}
\end{figure}

Figure~\ref{fig:graphQmu4} displays the graphs of $Q_{c,1}(4,\nu,\alpha_0)$ as functions of $\nu$ where the same color code as in Figures~\ref{fig:graphQnu1} and~\ref{fig:graphQnu2} is applied. Unlike these figures where all the plotted functions are 
monotonically decreasing, Figure~\ref{fig:graphQmu4} shows that $Q_{c,1}(4,\nu,\alpha_0)$ decreases monotonically for $\alpha_0=0$ and 1, but exhibit increasing behavoir for $\alpha_0=10,\,100$ and 500. This change of behavior is related to the fact that for small diffusion exponent $\nu$ (i.e. $\nu\,\ll\,\mu\,=\,4$) the critical value of $Q$ is strongly reduced for larger values of $\alpha$. In contrast, when $\nu$ is close to $\mu=4$, the strategy consisting in the choice of larger values of $\alpha$ becomes limited. 

\subsection{Minimum $\alpha$ for survival: numerical estimation}

In this section, our goal is to determine the minimum value of $\alpha$ in the initial conditions~(\ref{rho0fam1}) and~(\ref{rho0fam2}) that ensures the persistence of a population whose evolution is governed by the exponents $\mu=\mu_0$ and $\nu=\nu_0$ when we assume that the habitat size is fixed and that a given total population is distributed within it.
 Thus, according to~(\ref{Q}) the value of $Q$ is fixed. Let us denote $\widehat{Q}_c(\mu_0,\nu_0)=Q_{c,1}(\mu_0,\nu_0,0)=Q_{c,2}(\mu_0,\nu_0,0)$. If $Q\geq \widehat{Q}_c(\mu_0,\nu_0)$ population survival is guaranteed, otherwise, the only adjustable parameter that can ensure survival is $\alpha$. As it is shown in the previous subsection, a more concentrated initial population (i.e., larger $\alpha$) increases the likelihood of persistence. Thus, for a given $Q_0< \widehat{Q}_c(\mu_0,\nu_0)$, we look for the minimum value of $\alpha$ that allows population survival. We denote this threshold as $\alpha_{\min,j}(\mu_0,\nu_0,Q_0)$ where $j=1$ for~(\ref{rho0fam1}) and $j=2$ for~(\ref{rho0fam2}).

Analogously, for simplicity of notation, we outline the procedure for estimating this cut-off value by considering the first family of initial conditions. Following the definition of $\alpha_{\min,1}(\mu_0,\nu_0,Q_0)$, we start with $\alpha^{(0)}$ large enough so that the initial condition~(\ref{rho0fam1}) with $\alpha=\alpha^{(0)}$ leads to a total population that, after a possible transient phase, becomes a monotonically increasing function that reaches arbitrarily large values (with possible blow-up). Next, we choose a step size $\Delta\,\alpha$ small enough compared with $\alpha^{(0)}$ and consider the sequence of values $\alpha^{(r)}=\alpha^{(r-1)}-\Delta\,{\alpha}$, for $r\in\mathbb{N}$. For increasing values of $r$, we apply the numerical algorithm to approximate the total population~(\ref{NTapprox}) corresponding to the initial condition~(\ref{rho0fam1}) with $\alpha=\alpha^{(r)}$. This process continues until, for a certain $R$, the approximated total population becomes an asymptotically decreasing function which describes the extinction of the population. We denote $\alpha^*=\alpha^{(R-1)}$ and $\alpha^{**}=\alpha^{(R)}$ and  approximate $\alpha_{\min,1}(\mu_0,\nu_0,Q_0)\approx\alpha^*$.  The same method is applied to the second family of initial conditions. 

\begin{figure}[H]
	\centering
        \includegraphics[width=14cm]{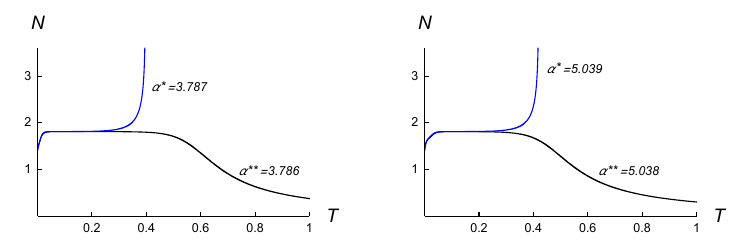}
        \caption{Evolution of the total population for equation~(\ref{prob2}) with $\mu_0=4$, $\nu_0=2$, $Q_0=2<\widehat{Q}_c(4,2)\approx4.467$ and initial condition~(\ref{rho0fam1}) (left) and~(\ref{rho0fam2}) (right). The left plot shows that the total population grows for $\alpha^*=3.787$ (blue line) while it decreases for $\alpha^{**}=3.786$ (black line), which leads to $\alpha_{\min,1}(4,2,2)\approx3.787$. Similarly, from the right plot it is observed how the total population grows for $\alpha^*=5.039$ (blue line) and decreases for $\alpha^{**}=5.038$ (black line), leading to $\alpha_{\min,2}(4,2,2)\approx5.039$.
         \label{fig:figalphamin}}
        \end{figure}

We illustrate the method for $\mu_0=4$, $\nu_0=2$ and  $Q_0=2$. Since $\widehat{Q}_c(4,2)\approx4.467$, it follows that $Q_0<\widehat{Q}_c(4,2)$.
The choice $\Delta\,\alpha=0.001$ results in $\alpha^*=3.787$ and $\alpha^{**}=3.786$ for the first family, and leads to $\alpha^*=5.039$ and $\alpha^{**}=5.038$ for the second one. Thus, we approximate  $\alpha_{\min,1}(4,2,2)\approx 3.787$ and $\alpha_{\min,2}(4,2,2)\approx 5.039$. Figure~\ref{fig:figalphamin} displays the approximate total population~(\ref{NTapprox}) for $\alpha^*$ (blue lines) and $\alpha^{**}$ (black lines). The left plot corresponds to the first family of initial conditions, whereas the right plot corresponds to the second family. 
        
\begin{figure}[H]
\centering
\includegraphics[width=12cm]{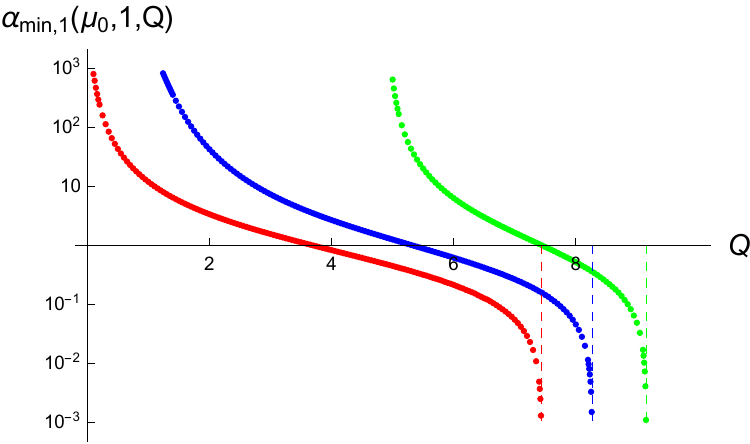}
\caption{Plots of $\alpha_{\min,1}(\mu_0,1,Q)$ for $\mu_0=2$ (green), $\mu_0=3$ (blue) and $\mu_0=4$ (red). Logarithmic scale is used in the vertical axis so that the vertical asymptotes indicate the critical values of $Q$: $\widehat{Q}_{c}(2,1)\approx9.153$ (green), $\widehat{Q}_{c}(3,1)\approx8.269$ (blue) and $\widehat{Q}_{c}(4,1)\approx7.443$ (red).}
\label{fig:graphalphamin1nu1}
\end{figure}

Figure~\ref{fig:graphalphamin1nu1} displays the graph of $\alpha_{\min,1}(\mu_0,1,Q)$ for $\mu_0=2$ (green), $\mu_0=3$ (blue) and $\mu_0=4$ (red). We point out that when the scheme to estimate $\alpha_{\min,1}(\mu_0,1,Q)$ is applied to increasing values of $Q$, the step size $\Delta\,{\alpha}$ may need to be reduced ensuring that it remains sufficiently small compared to the estimated value of 
$\alpha_{\min,1}(\mu_0,1,Q)$. 

The red curve, represented by a set of computed points, reveals that for the exponents $\mu_0=4$ and $\nu_0=1$, even for very small values of $Q$, it is possible to select appropriate values of $\alpha$ (and consequently adopt suitable strategies to concentrate the initial population) in order to ensure population survival. The minimum required value of $\alpha$ to prevent extinction decreases as $Q$ increases, tending to zero as $Q$ approaches to $\widehat{Q}_c(4,1)\approx7.443$. 

On the other hand, the blue curve indicates that for $\mu_0=3$ and $\nu_0=1$, the minimum $\alpha$ required for population persistence becomes very large for values of $Q$ closed to 1 (for example $\alpha_{\min}(3,1,1.3)\approx600.5$). Note, that very large $\alpha$ corresponds to an initial distribution that it is almost zero throughout the entire interval $[-\tfrac{1}{2},\tfrac{1}{2}]$ except in a very small subinterval around the origin, where the initial density reaches extremely high values. As a consequence, numerical simulations become less reliable. Moreover, it would not be realistic to assume that such an initial distribution could be practically implemented in the habitat. These observations suggest that a feasible strategy to prevent extinction can be applied for $Q\gtrapprox1.3$.
Again, $\alpha_{\min,1}(3,1,Q)$ decreases as $Q$ increases and tends to zero as $Q$ approaches $\widehat{Q}_c(3,1)\approx8.269$.

A similar interpretation applies to the green curve, where $\alpha_{\min,1}(2,1,5)\approx651.7$ and $\widehat{Q}_c(2,1)\approx9.153$. 

\begin{figure}[H]
\centering
\includegraphics[width=12cm]{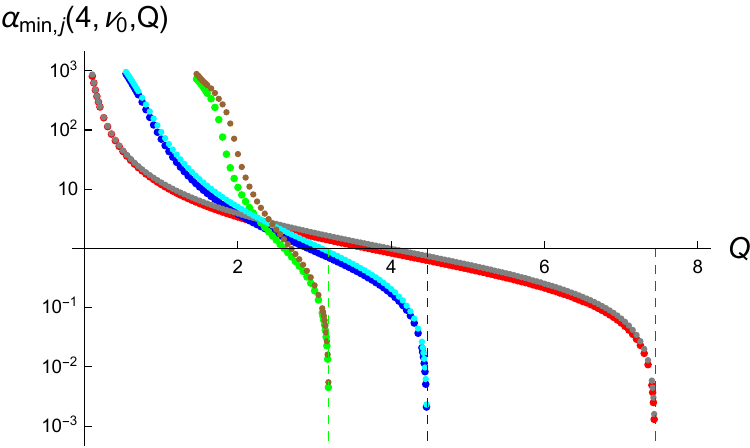}
\caption{Plots of $\alpha_{\min,j}(4,\nu_0,Q)$ for $j=1$ with $\nu_0=1$ (red), $\nu_0=2$ (blue), $\nu_0=3$ (green) and $j=2$ with $\nu_0=1$ (gray), $\nu_0=2$ (cyan), $\nu_0=3$ (brown). Logarithmic scale is used in the vertical axis so that the vertical asymptotes indicate the critical values of $Q$: $\widehat{Q}_c(4,1)\approx7.443$ (red), $\widehat{Q}_c(4,2)\approx4.467$ (blue) and $\widehat{Q}_c(4,3)\approx3.185$ (green).}
\label{fig:graphalphamin12mu4}
\end{figure}

\begin{figure}[H]
\centering
\includegraphics[width=14cm]{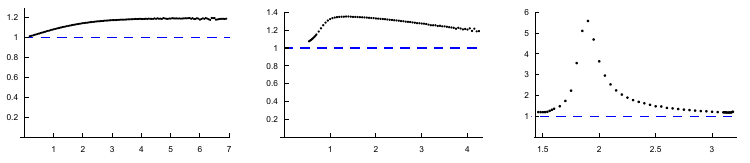}
\caption{Plots of the quotient $\alpha_{\min,2}(4,\nu_0,Q)/\alpha_{\min,1}(4,\nu_0,Q)$ for $\nu_0=1$ (left), $\nu_0=2$ (middle) and
$\nu_0=3$ (right).}
\label{fig:graphquotiensalpha}
\end{figure}

Figure~\ref{fig:graphalphamin12mu4} displays the graphs of $\alpha_{\min,1}(4,\nu_0,Q)$ for $\nu_0=1$ (red), $\nu_0=2$ (blue) and $\nu_0=3$ (green) as well as the graphs of $\alpha_{\min,2}(4,\nu_0,Q)$ for $\nu_0=1$ (gray), $\nu_0=2$ (cyan), $\nu_0=3$ (brown). Note that the red curve is the same as in Figure~\ref{fig:graphalphamin1nu1}. Each one of the dotted curves can be interpreted in the same way as in Figure~\ref{fig:graphalphamin1nu1}. Moreover, Figure \ref{fig:graphalphamin12mu4} shows that $\alpha_{\min,2}(4,\nu_0,Q)$ is slightly larger than $\alpha_{\min,1}(4,\nu_0,Q)$ for $\nu_0=1$ and $\nu_0=2$. This means that when the initial distribution is concentrated around $X=1/6$ a slightly larger $\alpha$ (i.e. a slightly more concentrated initial distribution) is required to ensure population survival. In contrast, for $\nu_0=3$ the minimum $\alpha$ necessary to guarantee population persistence under the initial distribution~(\ref{rho0fam2}) can be significantly larger than the one required  when we use the initial distribution~(\ref{rho0fam1}). These results are further illustrated in Figure~\ref{fig:graphquotiensalpha} where the quotient
\begin{equation}\label{quot}
\frac{\alpha_{\min,2}(4,\nu_0,Q)}{\alpha_{\min,1}(4,\nu_0,Q)},
\end{equation}
is displayed for $\nu_0=1$ (left), $\nu_0=2$ (middle) and $\nu_0=3$ (right). The figure shows that the quotient remains below 1.2 for $\nu_0=1$ and below 1.4 for $\nu_0=2$, whereas, it can reach values close to 5 for $\nu_0=3$.

Finally, both Figures~\ref{fig:graphalphamin1nu1} and~\ref{fig:graphalphamin12mu4} reveal that the choice of an appropriate $\alpha$ to ensure population survival for $Q$ below $\widehat{Q}_c(\mu,\nu)$ becomes less feasible as $\mu$ and $\nu$ become closer (with $\mu>\nu$). This fact is in agreement with equation~(\ref{Qc}) which states that
\begin{equation}\label{Qcalpha}
Q_{c,j}(\mu,\mu,\alpha)=\frac{\pi^2}{\mu},\quad j=1,2,\quad \mu>0,
\end{equation}
and consequently the critical value does not depend on $\alpha$.

\subsection{Critical habitat size and critical total population: numerical estimations \label{sec:lcn0c}}

In Section~\ref{sec:critval} we have discussed the condition for population survival in terms of the parameter $Q$, defined in~(\ref{Q}), which from the mathematical point of view, is the natural parameter to approach the problem. However, from the ecological point of view, one would be mainly interested in the conditions on either the habitat size (for a given total population) or the total population (for a given habitat size). As we have previously indicated in Section~\ref{sec:model}, the condition on the habitat size applies for $\mu$ and $\nu$ such that $\mu\,\neq\,\nu\,+\,2$, and for our families of initial distributions takes the form
\begin{equation}\label{hscon}\begin{array}{llll}
l\,\geq\,l_{c,j}(\mu,\nu,\alpha)&\;\mbox{ for }\;&\mu\,<\,\nu\,+\,2,&\;j\,=\,1,\,2,\\
l\,\leq\,l_{c,j}(\mu,\nu,\alpha)&\;\mbox{ for }\;&\mu\,>\,\nu\,+\,2,&\;j\,=\,1,\,2,
\end{array}\end{equation}
with $l_{c,j}(\mu,\nu,\alpha)$ being given by
\begin{equation}\label{lcfam}
l_{c,j}(\mu,\nu,\alpha)\,=\,\left(\frac{D}{a}Q_{c,j}(\mu,\nu,\alpha)\right)^{\frac{1}{-\mu+\nu+2}}n_0^{\frac{\mu-\nu}{\mu-\nu-2}}.
\end{equation}
Analogously, it is immediate to derive from~(\ref{Q}) and~(\ref{fcond}) that the survival condition on the total population for $\mu\,>\,\nu$ and a fixed habitat size $l$ is given by
\begin{equation}\label{n0cfam}
n_0\,\geq\,n_{0,c,j}(\mu,\nu,\alpha)=\left(\frac{D}{a}Q_{c,j}(\mu,\nu,\alpha)\right)^{\frac{1}{\mu-\nu}}l^{\frac{\mu-\nu-2}{\mu-\nu}},\quad j\,=\,1,\,2.
\end{equation}
Note that the dependences of $l_{c,j}(\mu,\nu,\alpha)$ and $n_{0,c,j}(\mu,\nu,\alpha)$ on the exponents $\mu$ and $\nu$ are considerably involved.

\begin{figure}[H]
\centering
\includegraphics[width=6cm]{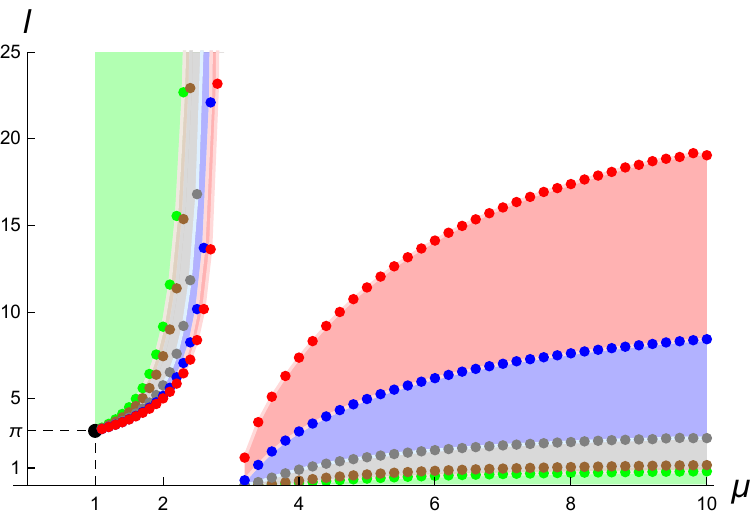}\hspace{5mm}
\includegraphics[width=6cm]{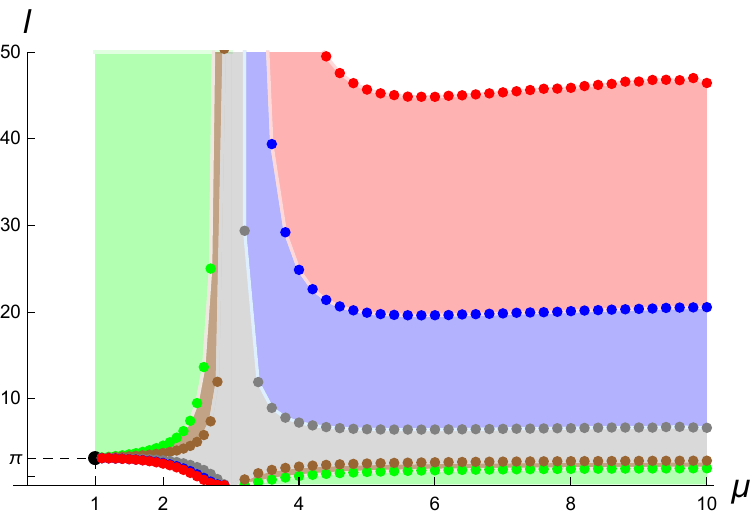}
\caption{Critical habitat size $l_{c,1}(\mu,1,\alpha_0)$ for $n_0=1$ (left) and $n_0=2$ (right) as functions of $\mu$ for $\alpha_0=0$ (green), 1 (brown), 10 (gray), 100 (blue) and 500 (red), when $D=a=1$.  The green shaded region corresponds to the survival region for $\alpha_0=0$ (i.e. the homogeneous initial distribution). As $\alpha_0$ increases, the survival region enlarges, adding progressively  the brown shaded region for $\alpha_0=1$, the gray shaded region for $\alpha_0=10$, the blue shaded region for $\alpha_0=100$ and the red shaded region for $\alpha_0=500$. The black point $(1,\pi)$ corresponds to $l_{c,1}(1,1,\alpha)=\pi$ for all $\alpha$ and $n_0$, when $D=a=1$.}
\label{fig:graphlcnu1}
\end{figure}

Figure~\ref{fig:graphlcnu1} displays the graphs of $l_{c,1}(\mu,1,\alpha_0)$ for $n_0=1$ (left) and $n_0=2$ (right) as functions of $\mu$ for $\alpha_0=0$ (green), $\alpha_0=1$ (brown), $\alpha_0=10$ (gray), $\alpha_0=100$ (blue) and $\alpha_0=500$ (red), when $D=a=1$. The shaded regions correspond to the survival regions in the parameter plane $(\mu,l)$. More precisely, the green shaded region is the survival region for the homogeneous initial distribution ($\alpha_0=0$) and as $\alpha_0$ increases the survival region enlarges progressively, adding the brown shaded region for $\alpha_0=1$, the gray shaded region for $\alpha_0=10$, the blue shaded region for $\alpha_0=100$ and the red shaded region for $\alpha_0=500$. It is worth noting that the green shaded region in Figure~\ref{fig:graphlcnu1} (left) corresponds to Figure $5$b in \cite{colombo2018nonlinear}. Furthermore, one can observe that, for sufficiently large values of $\mu$ (approximately $\mu\gtrapprox4$ and therefore $\mu>\nu+2$), a higher value of $\alpha_0$ considerably increases the critical value $l_{c,1}(\mu,1,\alpha_0)$ and consequently facilitates population survival. However, if $\mu<\nu+2$ or $\mu\gtrapprox\nu+2$ a considerably larger $\alpha$ is required to achieve a noticeable reduction in $l_c(\mu,1,\alpha_0)$. When $n_0=2$, the situation changes drastically, increasing the survival regions for large values of $\alpha_0$ in a deleted neighbourhood of $\mu=\nu+2=3$. This occurs, on the one hand, because for $\mu<3$, $l_{c,1}(\mu,1,\alpha_0)$ is a decreasing function that tends to zero as $\mu$ approaches $3$ when $\alpha_0=10$, $\alpha_0=100$ and $\alpha_0=500$; while for $\mu$ slightly greater than $3$, $l_{c,1}(\mu,1,\alpha_0)$ takes considerably large values when $\alpha_0=10$, $\alpha_0=100$ and $\alpha_0=500$.

\begin{figure}[H]
\centering
\includegraphics[width=8cm]{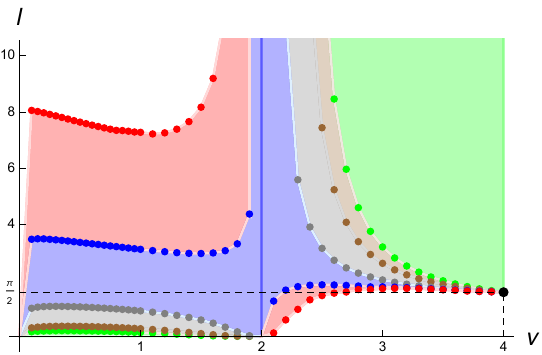}
\caption{Critical habitat size $l_{c,1}(4,\nu,\alpha_0)$ for $n_0=1$ as functions of $\nu$ for $\alpha_0=0$ (green), 1 (brown), 10 (gray), 100 (blue) and 500 (red), when $D=a=1$.  The green shaded region corresponds to the survival region for $\alpha_0=0$ (i.e. the homogeneous initial distribution). As $\alpha_0$ increases, the survival region enlarges, adding progressively  the brown shaded region for $\alpha_0=1$, the gray shaded region for $\alpha_0=10$, the blue shaded region for $\alpha_0=100$ and the red shaded region for $\alpha_0=500$. The black point $(4,\frac{\pi}{2})$ corresponds to $l_{c,1}(4,4,\alpha)=\pi$ for all $\alpha$ and $n_0$, when $D=a=1$.}
\label{fig:graphlcmu4}
\end{figure}

Figure~\ref{fig:graphlcmu4} displays the graphs of $l_{c,1}(4,\nu,\alpha_0)$ for $n_0=1$ as functions of $\nu$ using the same color code as in Figure~\ref{fig:graphlcnu1}, when $D=a=1$. As before, it should be noted that the green shaded region to the right in Figure~\ref{fig:graphlcmu4} could be compared to Figure $5$a in \cite{colombo2018nonlinear}, nevertheless the authors restrict their anaylisis to $\mu \leq 2$ and consequently, the behaviour displayed on the left of Figure~\ref{fig:graphlcmu4} can not be illustrated. For $n_0=1$, the behaviour of $l_{c,1}(4,\nu,\alpha_0)$ compared to $l_{c,1}(\mu,1,\alpha_0)$ differs significantly. This is due to the fact that $Q_{c,1}(4,\nu,\alpha_0)$ takes values below $1$ when $\alpha_0=100$ and $\alpha_0=500$ for values of $\nu$ arbitrarily close to $2$ (see Figure \ref{fig:graphQmu4} (right)). Consequently, $l_c(4,\nu,\alpha_0)$ attains large values when $\mu \lessapprox \nu+2$ and small values when $\mu \gtrapprox \nu+2$, which in turn results in a considerably expansion of the survival regions for large values of $\alpha_0$ when $\mu$ approaches $\nu+2$ from either side. In contrast, when $\nu$ is close to $\mu=4$, the strategy of considering large values of $\alpha_0$ becomes limited.

\begin{figure}[H]
\centering
\includegraphics[width=6cm]{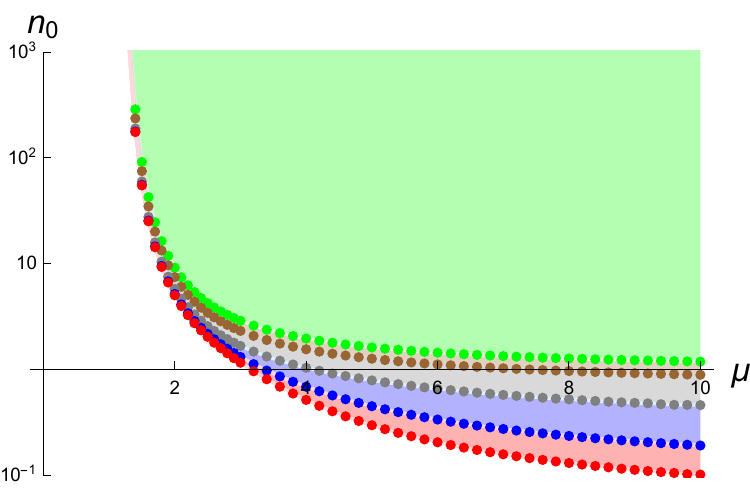}\hspace{5mm}
\includegraphics[width=6cm]{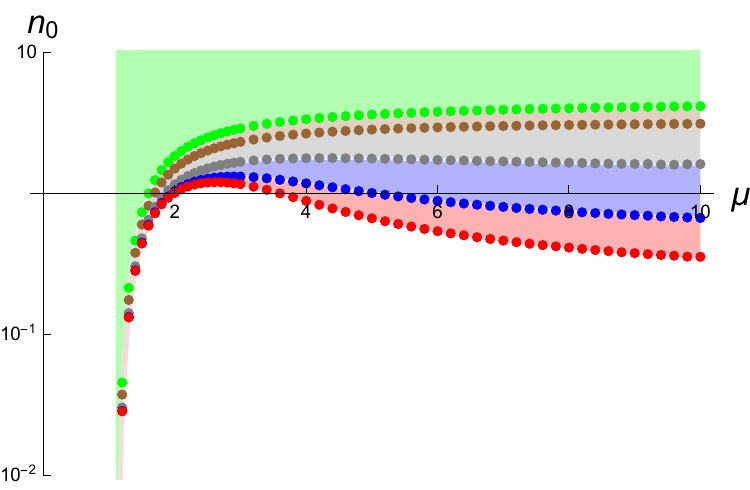}\\
\includegraphics[width=6cm]{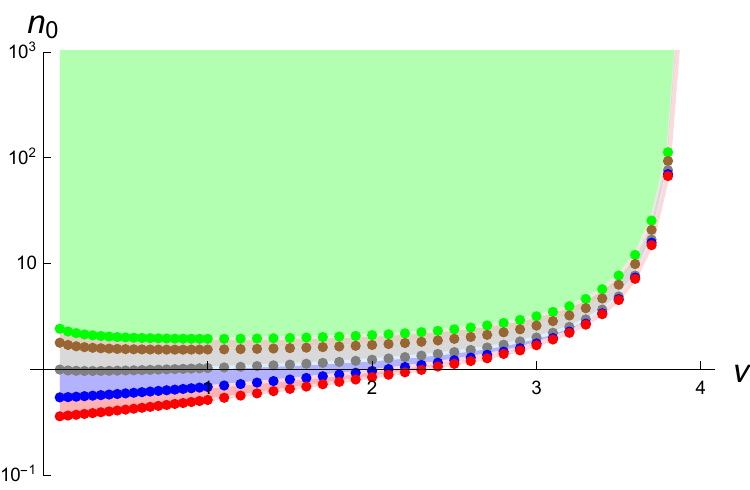}\hspace{5mm}
\includegraphics[width=6cm]{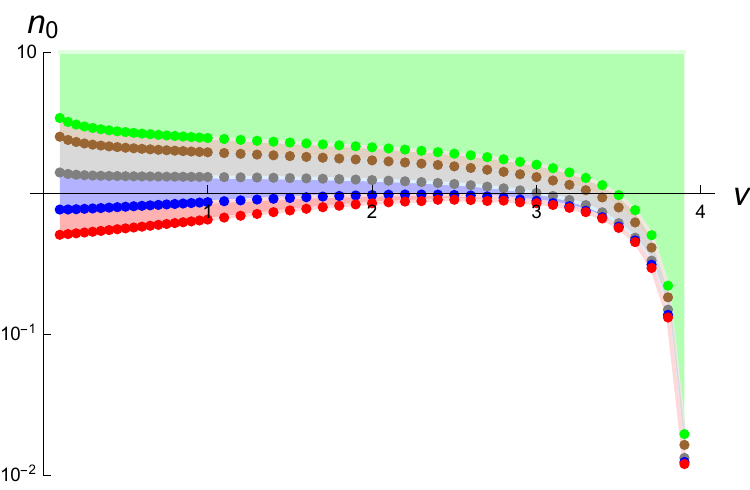}
\caption{Critical total initial population $n_{0,c,1}(\mu,1,\alpha_0)$ (top) for $l=1$ (left), $l=5$ (right)
and $n_{0,c,1}(4,\nu,\alpha_0)$ (bottom) for $l=1$ (left) and $l=2$ (right),  when $D=a=1$.  The green shaded region corresponds to the survival region for $\alpha_0=0$ (i.e. the homogeneous initial distribution). As $\alpha_0$ increases, the survival region enlarges, adding progressively the brown shaded region for $\alpha_0=1$, the gray shaded region for $\alpha_0=10$, the blue shaded region for $\alpha_0=100$ and the red shaded region for $\alpha_0=500$. A logarithmic scale is used on the vertical axis in the plots to highlight the differences among the different values of $\alpha_0$.}
\label{fig:graphn0cnu1}
\end{figure}

Figure~\ref{fig:graphn0cnu1} displays the graphs of $n_{0,c,1}(\mu,1,\alpha_0)$ (top) for $l=1$ (left) and $l=5$ (right) as functions of $\mu$, and the graphs of $n_{0,c,1}(4,\nu,\alpha_0)$ (bottom) for $l=1$ (left) and $l=2$ (right) as functions of $\nu$, when $D=a=1$, using the same color code as in Figure~\ref{fig:graphlcnu1} and with the shaded regions corresponding to the survival regions in the parameter plane $(\mu,n_0)$ (top) and the plane $(\nu,n_0)$ (bottom). A logarithmic scale is used on the vertical axis to highlight the differences among the different values of $\alpha_0$. Analogously, it can be observed that when the difference between $\mu $ and $\nu$ is large (approximately $\mu  \gtrapprox \nu+ 3$), an increasing value of $\alpha_0$ significantly reduces the critical values $n_{0,c,1}(\mu,1,\alpha_0)$ and $n_{0,c,1}(4,\nu,\alpha_0)$, and consequently facilitates population survival. Nevertheless, when $\mu \gtrapprox \nu$ a considerably larger $\alpha$ is required to achieve a noticeable reduction in $n_{0,c,1}(\mu,1,\alpha_0)$ and $n_{0,c,1}(4,\nu,\alpha_0)$, illustrating again that the strategy of using highly concentrated initial distributions becomes less effective. The variation observed between the graphs on the left and right in Figure~\ref{fig:graphn0cnu1} is attributed to the factor $l^{\frac{\mu-\nu-2}{\mu-\nu}}$ appearing in expression (\ref{n0cfam}).

\section{Concluding remarks \label{sec:con}}

We have analyzed the role of initial conditions for the persistence of a population in a bounded habitat $\left[-l/2,l/2\right]$ which evolves according to model (\ref{prob1})-(\ref{n0}). For this purpose, two families of initial conditions have been considered, one symmetric and centered at the origin (\ref{u01}), and the other centered at $\tfrac{l}{6}$ (\ref{u02}). 

A significant aspect of this work is the introduction of the parameter $Q$ (\ref{Q}), which not only allows a natural formulation of this problem but also offers a clearer interpretation of the underlying results. In particular, we have derived a survival condition in terms of $Q$ (\ref{fcond}) depending on a critical value $Q_{c}$ given by (\ref{Qc}), which holds universally for problem (\ref{prob1})-(\ref{n0}) under conditional persistence ($\mu \geq \nu$).

Numerical results have been presented to analyze the dependence of the critical values $Q_{c,j}$ on the exponents $\mu$, $\nu$ in the model, and the parameter $\alpha$ in the initial conditions (\ref{u01}) and (\ref{u02}). The results computed for both families of initial conditions lead us to conclude that initial condition (\ref{u01}) is slightly more propitious to ensure population survival (see Figure \ref{fig:graphdifQc}).

More importantly, we have developed a strategy to guarantee the persistence of a population whose evolution is determined by exponents $\mu=\mu_0$ and $\nu=\nu_0$, when it is assumed that the habitat size is fixed and a given total population is distributed within it. In particular, for a given $Q_0<\widehat{Q}_c$, we have determined the minimum value of $\alpha$ which ensures population survival (see Figures \ref{fig:graphalphamin1nu1} and \ref{fig:graphalphamin12mu4}). Numerical results show that when the exponents $\mu$ and $\nu$ differ significantly, population survival is easily achieved even for small values of $Q$. In contrast, when the exponents $\mu$ and $\nu$ are relatively close, the minimum $\alpha$ required for population persistence becomes very large for small values of $Q$, resulting in an impractical initial distribution which is almost zero except within a very small subinterval centered at the origin (\ref{u01}) or at $l/6$ (\ref{u02}). This agrees with the fact that for $\mu=\nu$ the critical value does not depend on $\alpha$.

Furthermore, by using the relationship between the habitat size (for a given total population) with respect to $Q$, we have determined the critical habitat size $l_{c,j}$ (\ref{lcfam}) required to ensure the survival of the population. More precisely, given a total population $n_0$ and $\mu \in (\nu,\nu+2)$, we have shown that the persistence of the population is guaranteed when the habitat size $l$ is bigger than the critial size $l_{c,j}$, whereas if $\mu>\nu+2$, the persistence of the population is guaranteed when the habitat size $l$ is smaller than $l_{c,j}$. It is worth noting that equations (\ref{Q}), (\ref{fcond}) can be applied to the particular case $\mu=\nu+2$, which was not previously considered in \cite{colombo2018nonlinear}. For this case, we have shown that equations (\ref{Q}), (\ref{fcond}) only impose a restriction (\ref{critcase}) on the total initial population. Notably, the introduction of the parameter $Q$ also reveals a significantly distinct behaviour in the critical habitat size $l_{c,j}$ when the diffusion term is nonlinear and with respect to the total initial population $n_0$ (see Figures \ref{fig:graphlcnu1} and \ref{fig:graphlcmu4}).

Analogously, by using the relationship between the total initial population (for a given habitat size) with respect to $Q_{c,j}$, we have determined the minimum total initial population $n_{0,c,j}$ (\ref{n0cfam}) that guarantees the persistence of the population.  Similarly, Figure \ref{fig:graphn0cnu1} shows the different behaviour of $n_{0,c,j}$ with respect to $\mu$, $\nu$ and $l$. Moreover, it can be observed again that when the exponents $\mu$ and $\nu$ differ significantly, an increasing value of $\alpha_0$ substantially reduces $n_{0,c,j}$, thereby facilitating the survival of the population. Nonetheless, when $\mu$ and $\nu$ are relatively close achieving a noticeable reduction in $n_{0,c,j}$ requires a considerably larger value of $\alpha_0$.



\section*{Acknowledgments}

This work was partially supported by grant PID2022-140451OA-I00 funded by Ministerio de Ciencia e Innovación/Agencia Estatal de investigación \linebreak (doi:10.13039/501100011033) and by ``ERDF A way of making Europe''.

\bibliographystyle{elsart-num-sort}
\bibliography{population_vf}
\end{document}